\documentclass[final,1p,times]{elsarticle}
\usepackage{subfigure}

\usepackage{amssymb,amsmath,amsthm}

\newtheorem{lemma}{Lemma}
\newtheorem{theorem}[lemma]{Theorem}

\newtheorem{example}[lemma]{Example}
\newtheorem{conjecture}{Conjecture}

\begin{document}

\begin{frontmatter}



\title{MIRROR-CURVES AND KNOT MOSAICS}


\author{Slavik Jablan}

\address{The Mathematical Institute, Knez Mihailova
36, P.O.Box 367, 11001 Belgrade, Serbia, E-mail: sjablan@gmail.com}

\author{Ljiljana Radovi\'c}

\address{University of Ni\v s, Faculty of
Mechanical Engineering, A.~Medvedeva 14, 18 000 Ni\v s, Serbia,
E-mail: ljradovic@gmail.com}

\author{Radmila Sazdanovi\'c}

\address{Department of Mathematics, University of
Pennsylvania, 209 South 33rd Street, Philadelphia, PA 19104-6395,
USA, E-mail: radmilas@gmail.com}

\author{Ana Zekovi\' c}

\address{Zeta System, Golsvortijeva 1, 11000 Belgrade, Serbia, E-mail: ana@zeta.rs}

\begin{abstract}
Inspired by the paper on quantum knots and knot mosaics \cite{23}
and grid diagrams (or arc presentations), used extensively in the
computations of  Heegaard-Floer knot homology \cite{2,3,7,24}, we
construct the more concise representation of knot mosaics and grid diagrams via
mirror-curves. Tame knot theory is equivalent to knot mosaics \cite{23}, mirror-curves, and grid
diagrams \cite{3,7,22,24}. Hence, we introduce codes for mirror-curves treated as knot or link diagrams
placed in rectangular square grids, suitable for software implementation. We
provide tables of minimal mirror-curve codes for knots and links obtained from rectangular grids of size $3\times 3$ and $p\times 2$
($p\le 4$), and describe an efficient algorithm for  computing  the Kauffman
bracket and L-polynomials \cite{18,19,20} directly from mirror-curve representations.

\end{abstract}


\end{frontmatter}


\section{Introduction}
\label{}

Mirror-curves originated from matting, plaiting, and basketry.
They appear in arts of different cultures (as Celtic knots,
Tamil threshold designs, Sona sand drawings...), as well as in works of Leonardo
and D\"urer \cite{1,4,5,13,14,15,16,18}. P.~Gerdes recognized their
deep connection with the mathematical algorithmic-based structures:
knot mosaics, Lunda matrices, self-avoiding curves, and cell-automata
\cite{13,14,15,16}. Combinatorial complexity of Sona sand drawings
is analyzed by M.~Damian {\it et all} \cite{9} and E.D.~Demaine {\it et all} \cite{10}.

Mirror-curves are constructed out of rectangular square grids, denoted by $RG[p,q]$,
of dimensions $p$, $q$ ($p,q\in N$). First we connect the midpoints of adjacent edges of
$RG[p,q]$ to obtain a $4$-valent graph: every vertex of this graph is
incident to four edges, called {\it steps}. Next, choose a starting point and traverse
the curve so that we leave each vertex via the middle outgoing edge.
Returning to the starting point, is equivalent to closing a path called a {\it component}.
If we return to the starting point without traversing all of the steps, we choose
a different one and repeat the process until every step is used exactly once.
A {\it mirror-curve} in $RG[p,q]$ grid is the set of all components.
To obtain a knot or a link diagram from a mirror-curve we introduce the  ``over-under'' relation, turning each
vertex to the crossing, i.e., we choose a pair of collinear steps (out of two) meeting at a vertex to be the overpass \cite{18,19,20,25}.

Mirror-curves can also be obtained from the following physical model which, in a way, justifies
their name: assume that the sides of our rectangular square grid $RG[p,q]$ are made
of mirrors, and that additional internal two-sided mirrors are placed between the square
cells, coinciding with an edge, or perpendicular to it in its midpoint. If a ray of light is emitted
from one edge-midpoint at an angle of $45^\circ ,$ it will eventually
come back to its starting point, closing a component after series of reflections. If some steps remained
untraced, repeat the whole procedure starting from a different point.

Through the rest of the paper the term ``mirror-curves'' will be used for labeled mirror-curves.
Hence, all crossings will be signed, where $+1$ corresponds to the positive, and $-1$ to negative crossings.

\begin{theorem} $[15]$
The number of components of a knot or link $L$ obtained from a rectangular grid $RG[p,q]$ without
internal mirrors is $c(L)=GCD(p,q).$
\end{theorem}

The web-Mathematica computations with mirror-curves are available at
the address

\medskip

{\tt http://math.ict.edu.rs:8080/webMathematica/mirror/cont.htm}

\section{Coding of mirror-curves}
\label{}

\begin{figure}[th]
\centering
{ \includegraphics[scale=0.6]{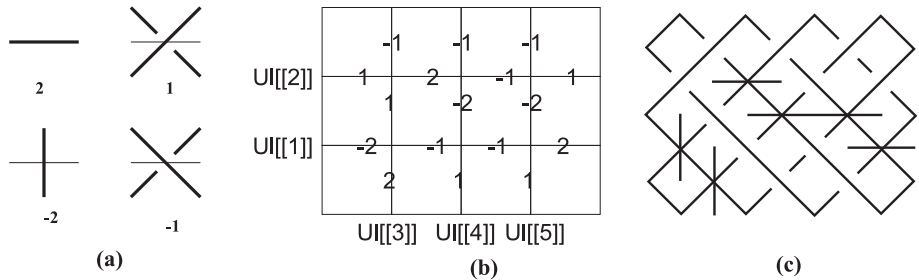}}
\caption{(a) Edge labeling; (b) labeled $RG[3,2]$; (c) the mirror-curve corresponding to the code $Ul$.}\label{f1.1}
\end{figure}

Mirror-curve is constructed on a rectangular grid $RG[p,q]$ with
every internal edge labeled $1$,  $-1$, $2$, and $-2$, where
 $+1$ and $-1$ denote, respectively, a positive and negative crossing in the
middle point of the edge, see Figure~\ref{f1.1}a, while
 $2$ and $-2$ denote a two-sided mirror containing the middle point of an edge,
either collinear or perpendicular to it.
The code for the mirror-curves can be given in matrix form, containing labels
of internal edges corresponding to rows and columns of the $RG[p,q]$. For
example,  the code
$$Ul=\{\{-2,-1,-1,2\},\{1,2,-1,1\},\{2,1,-1\},\{1,-2,-1\},\{1,-2,-1\}\}.$$
corresponds to the mirror-curve on Figure~\ref{f1.1}c, based on the labeled rectangular grid $RG[3,2]$ shown in Figure~\ref{f1.1}b.

Our convention is the natural one: we list labels in the rows from left to right, and in the
columns from bottom to the top.

\section{Reduction of mirror-curves}

Labeled mirror-curves represent knot and link (shortly $KL$) diagrams. In this section we consider Reidemeister moves, expressed in the language of mirror-curves.

\begin{figure}[th]
\centering
{ \includegraphics[scale=0.4]{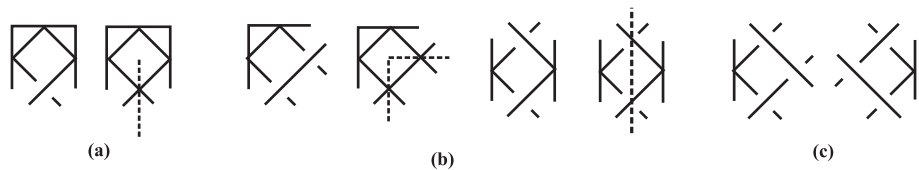}}
\caption{(a) Reidemeister move $R$I; (b) Reidemeister
move $R$II; (c) Reidemeister move $R$III, with additional mirrors in $R$I and
$R$II denoted by dotted lines.\label{f1.3}}
\end{figure}

The Reidemeister move $R$I is equivalent to replacing crossing by the mirror $-2$ (i.e., $\pm 1 \rightarrow
-2$), see Figure~\ref{f1.3}a.

Reidemeister move $R$II is the replacement of two neighboring crossings of the same sign by two perpendicular or collinear mirrors shown on Figure~\ref{f1.3}b, and Reidemeister move $R$III is illustrated in Figure~\ref{f1.3}c.

\begin{figure}[th]
\centering
{ \includegraphics[scale=0.5]{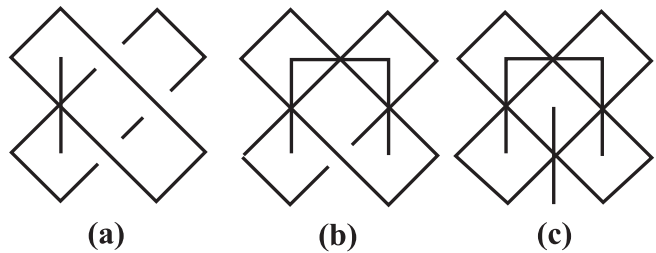}}
\caption{A sequence of Reidemeister moves reducing $3$-crossing diagram of an unknot to the minimal one.\label{f1.4}}
\end{figure}

Notice that every unknot or unlink can be reduced to the code containing only
labels $2$ and $-2$. For example, the non-minimal diagram of an unknot with three crossings on Figure  \ref{f1.4}a, given by the code $Ul = \{\{-2, -1\}, \{1, 1\}\}$, can be reduced using the second Reidemeister move $R$II
applied to the upper right crossings, to  $Ul=\{\{-2, -2\},
\{1, -2\}\}$ on Figure~\ref{f1.4}b. This code can be reduced further using the first Reidemeister move  $R$I applied to the remaining
crossing, yielding the minimal code of the unknot in $RG[2,2]$: $Ul = \{\{-2, -2\}, \{2, -2\}\}$.

Minimal diagrams of mirror-curves correspond to codes with the minimal number of $\pm 1$ labels. Minimal mirror-curve codes of alternating knots and links contain either $1$'s or $-1$'s, but not both of them.

\begin{figure}[th]
\centering
{ \includegraphics[scale=0.6]{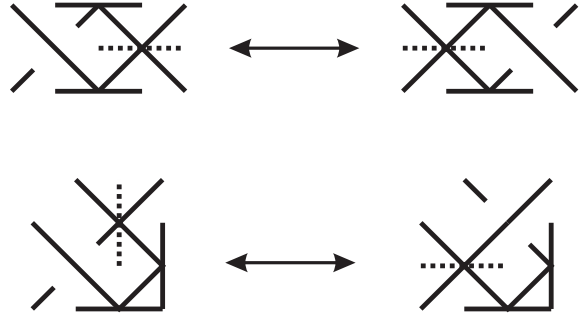}}
\caption{Mirror moves that can be useful for simplifying the reduction process.\label{f1.5}}
\end{figure}

Next we consider several examples to illustrate the reduction process.
Sometimes it is useful to use topological intuition to simplify the reduction,
such as the mirror-moves shown in Figure~\ref{f1.5}, where the repositioned mirror is shown by a dotted line.

K.~Reidemeister proved that any two different diagrams of the same knot or link are related by a finite sequence of Reidemeister moves, but there are no algorithms prescribing the order in which they can be used. Similarly, we have no algorithms for reducing mirror-curve codes. In particular,
we can not guarantee that we can obtain the minimal code without increasing the size of the rectangular grid.

\begin{example} This is the reduction sequence for the $2$-component link shown in Figure~\ref{f1.6}a, determined by the following code
$$\{\{-2,-1,-1,2\},\{1,2,-1,1\},\{2,1,-1\},\{1,-2,-1\},\{1,-2,-1\}\}$$
\noindent resulting in the unlink.
First we apply the first Reidemeister move $R$I to the right lower crossing in Figure~\ref{f1.6}b, and three moves
$R$II, in order to obtain the code
$$\{\{-2, -2, -1, 2\}, \{-2, 2, 2, -2\}, \{2, 1, -2\}, \{-2, -2,
  -1\}, \{-2, -2, -2\}\},$$
then the mirror-move to the first mirror in the upper row
and obtain the code corresponding to the Figure~\ref{f1.6}c:
 $$\{\{-2, -2, -1, 2\}, \{-2, 2, -1,-2\}, \{2, 1,-1\}, \{-2, -2, -2\}, \{-2, -2, -2\}\}.$$
Next we perform two Reidemeister moves $R$I to obtain Figure~\ref{f1.6}d, and the code
  $$\{\{-2, -2, -1, 2\}, \{-2, 2, 2, -2\}, \{2, 1, -1\}, \{-2, -2, -2\}, \{-2, -2, -2\}\},$$
and the link shown in Figure  \ref{f1.6}e:
 $$\{\{-2, -2, 2, 2\}, \{-2, 2, 2, -2\}, \{2, 1, -1\}, \{-2, -2, -2\}, \{-2, -2, -2\}\}.$$
 Finally, the second Reidemeister move  $R$II eliminates the remaining two crossings to
  obtain the minimal code see Figure~\ref{f1.6}f, $$\{\{-2, -2, 2, 2\}, \{-2, 2, 2, -2\}, \{2, 2,
  2\}, \{-2, -2, -2\}, \{-2, -2, -2\}\}.$$

\begin{figure}[th]
\centering
{ \includegraphics[scale=0.6]{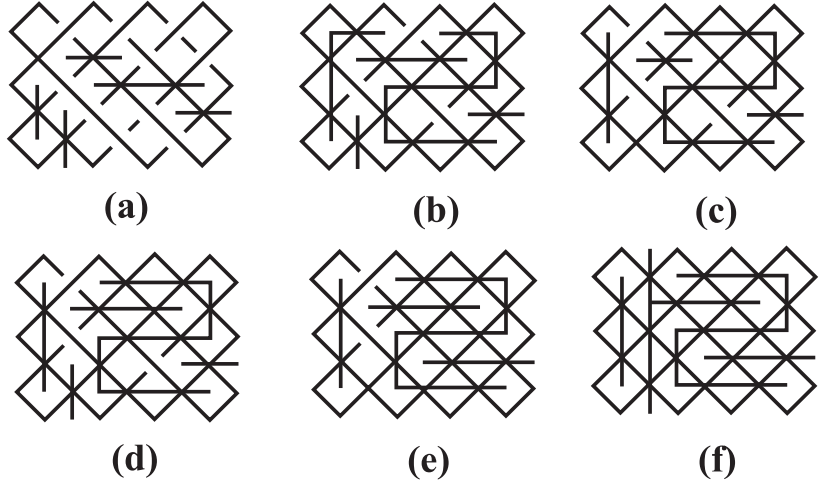}}
\caption{Reduction of two-component unlink
$Ul=\{\{-2,-1,-1,2\},\{1,2,-1,1\},$
$\{2,1,-1\},\{1,-2,-1\},\{1,-2,-1\}\}$.\label{f1.6}}
\end{figure}
\end{example}

Mirror-curve codes can be extended to virtual knots and links, by marking virtual crossings by zeros \cite{21}.

\section{Derivation of knots and links from mirror-curves}

Another interesting open problem is which knots and links
can be obtained from a rectangular grid $RG[p,q]$ of a fixed size.
To remove redundancies, we list each knot or link only once,
associated only with the smallest rectangular grid from which it
can be obtained.

Obviously, grid $RG[1,1]$ contains only the unknot, while from $RG[2,1]$ we can additionally
derive the trivial two-component unlink. In general, every rectangular grid $RG[p,1]$
contains the trivial $p$-component unlink.

In the rest of the paper, knots and links will be given by their classical notation and Conway
symbols \cite{6,18} from Rolfsen's tables \cite{25}. Links with more than 9 crossings
are given by Thistlethwaite's link notation \cite{3}.

\begin{figure}[th]
\centering
{ \includegraphics[scale=0.65]{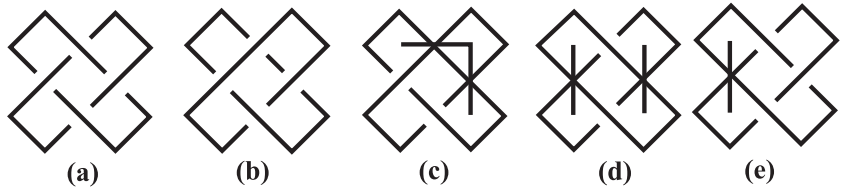}}
\caption{(a) Link $4$ ($4_1^2$); (b) non-minimal diagram of the Hopf
link $2$ ($2_1^2$); (c, d) two minimal diagrams of the Hopf link; (e)
minimal diagram of the trefoil knot $3$ ($3_1$). \label{f1.7}}
\end{figure}

Grid $RG[2,2]$ contains the following four knots and links shown in Figure~\ref{f1.7}: link
$4$ ($4_1^2$) given by the code $\{\{1, 1\}, \{1, 1\}\}$,
one non-minimal diagram of the Hopf link given by the code
$\{\{1, 1\}, \{1, -1\}\}$ which can be reduced to  the minimal diagram
$\{\{1, -2\}, \{1, -2\}\}$ using the second Reidemeister move $R$II, the symmetrical
minimal diagram of the Hopf link on Figure~\ref{f1.7}d, given by the code
$\{\{-2, -2\}, \{1, 1\}\}$, and the minimal diagram of trefoil (Figure~\ref{f1.7}e) given by the code $\{\{-2, 1\}, \{1, 1\}\}$.

Rectangular grid $RG[2,2]$ without internal mirrors, taken as the alternating
link, corresponds to the code which contains no $\pm 2$ and all $1$'s or exclusively $-1$'s. It
represents the link $4$ ($4_1^2$) (or its mirror image). Hence, the following two questions are equivalent:
\begin{itemize}
  \item  which
$KL$s can be obtained as mirror-curves from $RG[2,2];$
  \item which $KL$s can be obtained by substituting
crossings of the link $4$ ($4_1^2$) by elementary tangles $1$, $-1$,
$L_0$ and $L_\infty $, see Figure~\ref{f1.7}.
\end{itemize}

In analogy with the state sum model for the Kauffman bracket polynomial \cite{19}, where each crossing can be replaced by one of the two smoothings (resolutions) we can consider all possible states of a given rectangular grid $RG[p,q]$, corresponding to four different choices of placing a mirror $2$, $-2$, or one of the crossings $1$, $-1$ at the middle point of each edge. In this light, different mirror-curves obtained in this way can be thought of as all possible states of $RG[2,2]$, while the corresponding $KL$s can be viewed as all states of the link $4$ ($4_1^2$).

\begin{figure}[th]
\centering
{ \includegraphics[scale=0.6]{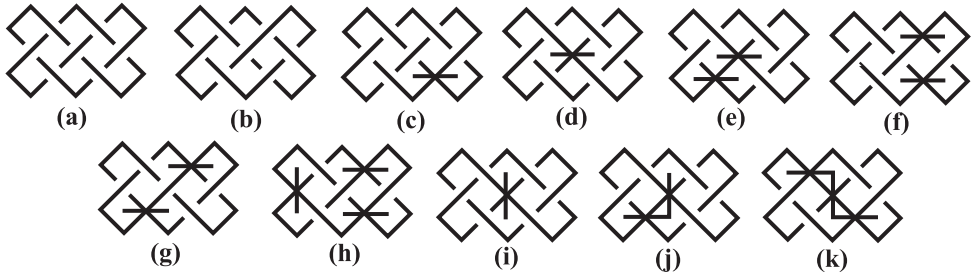}}
\caption{(a) Knot $3\,1\,3$ ($7_4$); (b) knot $4\,2$ ($6_1$); (c)
knot $3\,1\,2$ ($6_2$); (d) link 6 ($6_1^2$); (e) knot 5 ($5_1$);
(f) knot $3\,2$ ($5_2$); (g) Whitehead link $2\,1\,2$ ($5_1^2$); (h)
figure-eight knot $2\,2$ ($4_1$); (i) direct product of two trefoils
$3\#3$; (j) direct product of trefoil and Hopf link $3\#2$; (k)
direct product of two Hopf links $2\#2$.\label{f1.8}}
\end{figure}

From $RG[3,2]$ and its corresponding alternating knot $3\,1\,3$
($7_4$) given by the code $\{\{1, 1, 1\},$ $ \{1, 1\}, \{1, 1\}\}$
on Figure~\ref{f1.8}a, we obtain knots and links shown on Figure~\ref{f1.8}b--h:
\begin{center}
  \begin{tabular}{|l|l|}
  \hline
  $KL$ &Mirror-curve code\\
  \hline
    $4\,2$ & $\{\{1, 1, -1\}, \{1, 1\}, \{-1, -1\}\}$ \\
  $3\,1\,2$ ($6_2$) & $\{\{1, 1, 1\}, \{1, 1\}, $ $\{-2, 1\}\}$ \\
  $6$ ($6_1^2$)  &  $\{\{1, 2, 1\}, \{1,$ $ 1\}, \{1, 1\}\}$ \\
  $5$ ($5_1$) &  $\{\{1, 2, 1\}, \{-2,1\}, \{1, 1\}\}$\\
  $3\,2$ ($5_2$) & $\{\{1, 1, 1\},$ $ \{1, 1\},$ $ \{-2, -2\}\}$ \\
   $2\,1\,2$ ($5_1^2$) & $\{\{1, 1, 1\}, \{-2, 1\}, \{1, -2\}\}$  \\
    $2\,2$ ($4_1$) &  $\{\{-2, 1, 1\}, \{1, 1\}, \{-2, -2\}\}$\\
  \hline
\end{tabular}
\end{center}
\noindent and
the following composite knots and links shown on Figure~\ref{f1.8}i-k: direct product of two trefoils $3\#3$
given by the code $\{\{1, -2, 1\}, \{1, 1\},$ $ \{1, 1\}\}$, direct product of a trefoil and Hopf link $3\#2$ given by the
code $\{\{1, -2, 1\}, \{-2, 1\}, \{1, 1\}\}$, and direct
product of two Hopf links $2\#2$ given by the code $\{\{1, -2, 1\},
\{1,-2\},$ $\{-2,1\}\}$. In the case of composite knots and links we
can also obtain their non-alternating versions, e.g., $3\#(-3)$.

Alternating link $3\,1\,2\,1\,3$ ($L10a_{101}$ from Thistlethwaite's tables)
corresponds to $RG[4,2]$. The following prime knots and links can be obtained from
$RG[4,2]$: $5\,1\,3$ ($9_5$), $3\,1\,2\,1\,2$ ($9_{20}$),
$4\,1\,1\,3$ ($9_5^2$), $3\,1\,3\,2$ ($9_8^2$), $3\,1\,1\,1\,3$
($9_9^2$), $5\,1\,2$ ($8_2$), $4\,1\,3$ ($8_4$), $3\,1\,1\,1\,2$
($8_{13}$),  $8$ ($8_1^2$), $4\,2\,2$ ($8_3^2$), $3\,2\,3$
($8_4^2$), $3\,1\,2\,2$ ($8_5^2$), $2\,4\,2$ ($8_6^2$),
$2\,1\,2\,1\,2$ ($8_7^2$), $7$ ($7_1$), $5\,2$ ($7_2$), $2\,2\,1\,2$
($7_6$), $2\,1\,1\,1\,2$ ($7_7$),  $4\,1\,2$ ($7_1^2$), $3\,1\,1\,2$
($7_2^2$),  $2\,3\,2$ ($7_3^2$),  $2\,1\,1\,2$ ($6_3$), $3\,3$
($6_2^2$), and $2\,2\,2$ ($6_3^2$).

Moreover, we have a family of rational knots and links
corresponding to rectangular grids $RG[p,2]$ ($p\ge 3$), starting with $3\,1\,3$ ($7_4$), $3\,1\,2\,1\,3$
($L10a_{101}$), $3\,1\,2\,1\,2\,1\,3$, $\ldots $ given by their
minimal diagrams $3\,1\,3$, $(((1,(3,1),1),1),1,1,1)$,
$((1,(1,(1,(1,(1,3),1)),1)),1,1,1)$, $\ldots $ Rational knots, also known as $2$-bridge knots or $4$-plats\footnote{Knots or two component links obtained by a so-called horizontal closure of a braid on $4$ strings, with bottom connection points $A$, $B$, $C$, $D$, and the top connection points $A'$, $B'$, $C'$, $D'$, where we connect $A$ to $B$, $C$ to $D$, $A'$ to $B'$, and $C'$ to $D'$.}, form  the subset of
mirror-curves derived from rectangular grids $RG[p,2]$.

\begin{theorem}
All rational knots and links can be derived as mirror-curves from rectangular grids $RG[p,2]$ ($p\ge 2$).
\end{theorem}

\begin{figure}[th]
\centering
{ \includegraphics[scale=0.6]{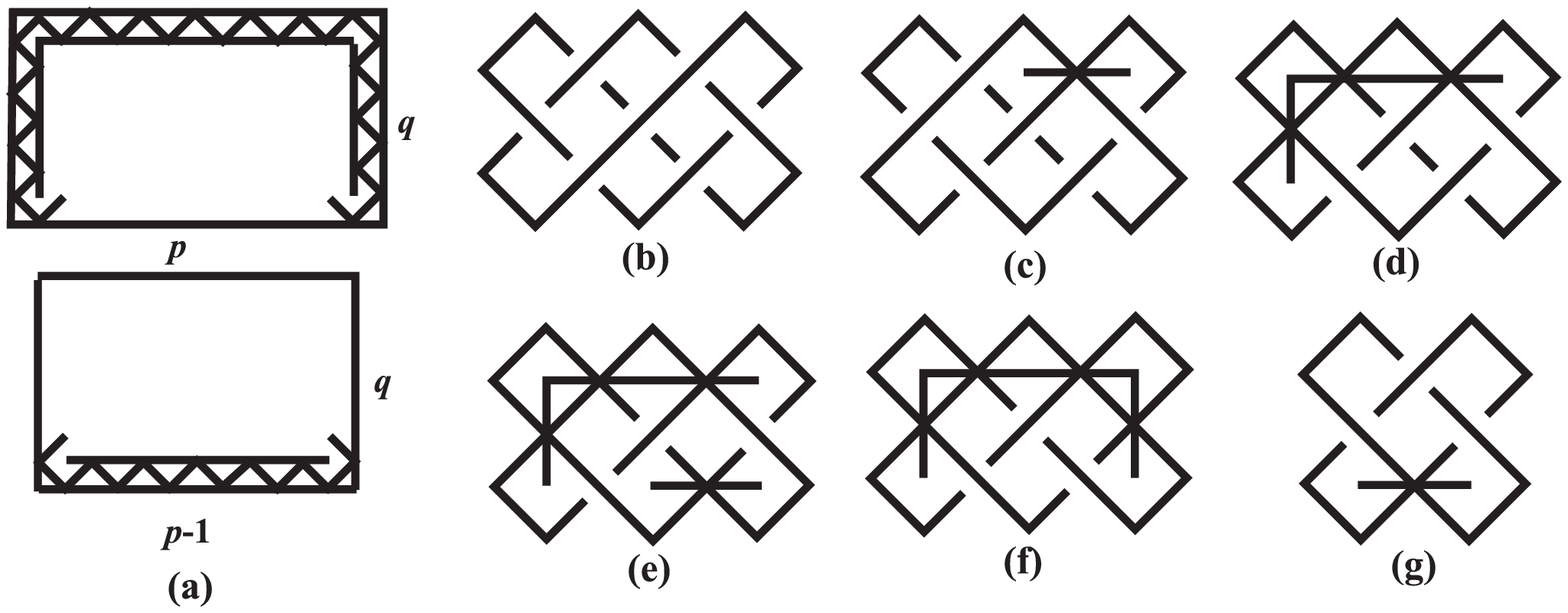}}
\caption{(a) Grid reduction by the all-over move; (b-g) six-step reduction of the knot $3\,-1\,3$ placed in the $RG[3,2]$ to the trefoil placed in its minimal grid $RG[2,2]$.}\label{f1.10}
\end{figure}

The reduction process we have described will not always result in the minimal rectangular grid for
representing a given $KL$ as a mirror-curve. Therefore we need a move that reduces the size of the grid,
so-called "all-over move", see Figure~\ref{f1.10}a, reducing the size of the grid from $RG[p,q]$ to $RG[p-1,q]$ while preserving the knot or link type.

The complete reduction of a non-minimal diagram of a trefoil, given
by the sequence of codes: $\{\{-1, 1, -1\}, \{-1, -1\}, \{-1,
-1\}\}$ $\rightarrow $ $\{\{1, 1, -1\}, \{1, -1\},$ $\{-1, -2\}\}$
$\rightarrow $ $\{\{-2, 1, -1\}, $ $ \{1, -2\},$ $\{-1, -2\}\}$
$\rightarrow $
 $\{\{-2, 1, -1\}, \{1, -2\}, \{-2,$ $-2\}\}$ $\rightarrow $
 $\{\{-2, 1, -2\}, \{1, -2\}, \{1, -2\}\}$ $\rightarrow $
 $\{\{-1, -1\}, \{-2, -1\}\}$,
including grid reduction from $RG[3,2]$ to $RG[2,2]$ in the last step, is illustrated in Figure~\ref{f1.10}b-g.

\begin{figure}[th]
\centering
{ \includegraphics[scale=0.65]{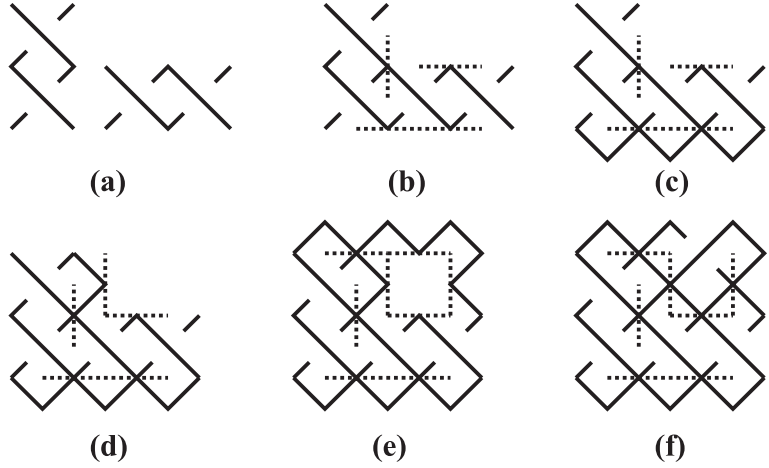}}
\caption{Construction of a mirror-curve diagram of the figure-eight knot
$4_1$ from its Conway symbol $2\,2$. \label{f1.11}}
\end{figure}

The next natural question is how to construct a mirror-curve representation of a knot
 or link given in Conway notation \cite{6,18,25}. We do not provide the general algorithm, but illustrate the process
 in the case of figure-eight knot $2\,2$. Knowing that the
figure-eight knot is obtained as a product of two tangles $2$, Figure~\ref{f1.11}a, we start by connecting two appropriate ends, see Figure~\ref{f1.11}b, and proceeding with completing the tangle $2\,2$ and its numerator closure.
In this process we are likely to obtain the empty regions, Figure~\ref{f1.11}e. They can be
incorporated in the construction by extending the mirror-curve across the empty region
included in our drawing by the Reidemeister move $R$I. This is achieved by deleting
a border mirror and changing the hole into a loop. Most often, mirror-curve representation
obtained in this way will not be the minimal one in terms of the grid size, so we
need to make further reductions\footnote{The simplest
way to obtain a mirror-curve from a given $KL$ is to use one of the programs {\it KnotAtlas}
\cite{3} or {\it gridlink} \cite{8} to construct a grid diagram of a
given link, then transform it into a mirror-curve, and make reduction at the end.}.

\begin{figure}[th]
\centering
{ \includegraphics[scale=0.35]{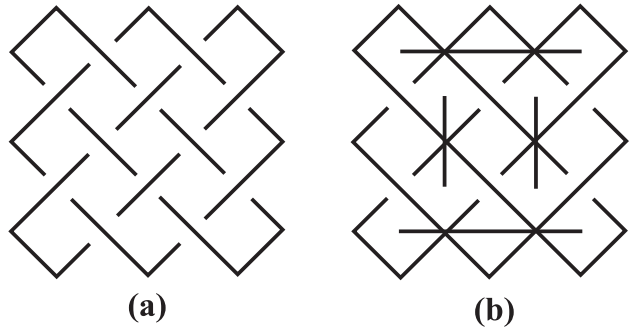}}
\caption{(a) $RG[3,3]$ with $3$-component link $8^*2:2:2:2$ (b)
non-alternating 3-component link $2,2,-2$ ($6_3^3$).\label{f1.12}}
\end{figure}

From $RG[3,3]$ and its corresponding alternating $3$-component link
$8^*2:2:2:2$ with $12$ crossings, given by the code $\{\{1, 1, 1\},
\{1, 1, 1\}, \{1, 1, 1\},$ $ \{1, 1, 1\}\}$, Figure~\ref{f1.12}a, we derive
many new knots and links, among them the smallest basic polyhedron -- Borromean
rings $6^*$ ($6_2^3$) given by the code $\{\{-1, -1, -2\},$ $ \{-2, -1,
-2\},$ $ \{-2, 2,$ $-2\}, \{-1,-1,-1\}\}$, see Figure~\ref{f1.13}b, and the first
non-alternating $3$-component link $2,2,-2$ ($6_3^3$) given by the
code $\{\{-1, -1, 1\},$ $ \{-1, -1, 1\}, \{-2, 2, -2\},$ $ \{-2, 2,
-2\}\}$ shown on Figure~\ref{f1.12}b.

Alternating link $8^*2:2:2:2$ corresponds to $RG[3,3]$, to which
we associate the following prime knots and links: $(2,2)\,(3\,1,-3\,1)$,
$(-5\,1,2)\,(2,2)$, $6^*-2.2.-2:4$, $6^*3.2.-3:2$, $6^*-3.-3\,
0::-3\, 0$, $2\,1\,2\,1\,1\,1\,2$ ($10_{44}$), $.4.2\,0$
($10_{85}$), $4\,1\,2\,1\,2$ ($L10a_{99}$), $.3:3\,0$
($L10a_{140}$), $6,2,2$ ($L10a_{145}$), $.2.3.2\,0$ ($L10a_{162}$),
$8^*2::2$ ($L10a_{163}$), $2\,0.2.2\,0.2\,0$ ($L10a_{164}$),
$(2\,1,-2\,1)\,(2,2)$ ($L10n_{73}$), $(3\,1,-2)\, (2,2)$
($L10n_{85}$), $(2,2)\,(4,-2)$ ($L10n_{86}$), $4,3\,1,-2$
($L10n_{92}$), $4,4,-2$ ($L10n_{93}$), $2\,0.-2.-2\,0.2\,0$
($L10n_{94}$), $3\,1,3\,1,-2$ ($L10n_{95}$), $4\,1\,2\,2$
($9_{11}$), $4\,1\,1\,1\,2$ ($9_{14}$), $2\,1\,3\,1\,2$ ($9_{17}$),
$2\,2\,1\,2\,2$ ($9_{23}$), $2\,1\,2\,1\,1\,2$ ($9_{27}$),
$2\,1\,1\,1\,1\,1\,2$ ($9_{31}$), $6\,1\, 2$ ($9_1^2$),
$2\,2\,1\,1\,1\,2$ ($9_{12}^2$), $5,2,2$ ($9_{13}^2$), $.4$
($9_{31}^2$), $.3.2\,0$ ($9_{35}^2$), $8^*2$ ($9_{42}^2$), $6\,2$
($8_1$), $3,3,2$ ($8_5$), $4\,1\,1\,2$ ($8_7$), $2\,3\,1\,2$
($8_8$), $2\,1,3,2$ ($8_{10}$) $2\,2\,2\,2$ ($8_{12}$),
$2\,2\,1\,1\,2$ ($8_{14}$), $.2.2\,0$ ($8_{16}$), $.2.2$ ($8_{17}$),
$8^*$ ($8_{18}$), $2\,1\,2\,1\,2$ ($8_7^2$), $2\,1\,1\,1\,1\,2$
($8_8^2$), $4,2,2$ ($8_1^3$), $3\,1,2,2$ ($8_2^3$), $(2,2)\,(2,2)$
($8_4^3$), $.3$ ($8_5^3$), $.2:2\,0$ ($8_6^3$), $4,2,-2$ ($8_7^3$),
$3\,1,2,-2$ ($8_8^3$), $(2,2)\,(2,-2)$ ($8_9^3$), $(2,2)\,-(2,2)$
($8_{10}^3$), $4\,3$ ($7_3$), $3\,2\,2$ ($7_5$), $2,2,2+$ ($7_1^3$),
$2\,3\,2$ ($7_3^2$), $3,2,2$ ($7_4^2$), $2\,1,2,2$ ($7_5^2$), $.2$
($7_6^2$), $2,2,2$ ($6_1^3$), $6^*$ ($6_2^3$), and $2,2,-2$
($6_3^3$).

\section{Knot mosaics, mirror-curves, grid diagram representations and tame knot theory}
\label{KM}
\begin{figure}[th]
\centering
{ \includegraphics[scale=0.65]{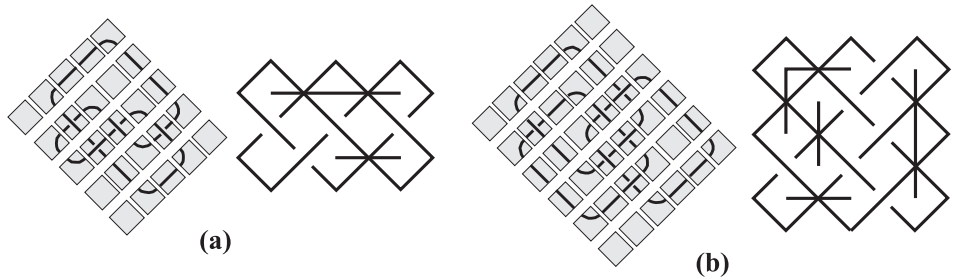}}
\caption{(a) Figure-eight knot and (b) Borromean rings from the
paper \cite{23} transformed into mirror-curves.\label{f1.13}}
\end{figure}

 Mirror-curves are equivalent to link mosaics: every link mosaic can be easily transformed
into a mirror-curve and {\it vice versa}. For example, the mosaics
of the figure-eight knot \cite{23} (pp.~$6$) and Borromean rings
\cite{23} (pp.~$7$) correspond to the mirror-curves on Figure~\ref{f1.13}
 and {\it vice versa}. Even more illustrative are knot mosaics
from the paper \cite{12} (pp.~$15$): first we rotate them by
$45^o$, cut out the empty parts, and add the two-sided mirrors in
appropriate places.

T.~Kuriya \cite{22} proved Lomonaco-Kauffman Conjecture
\cite{23}, showing that the tame knots are equivalent to knot mosaics, hence
also to mirror-curves. According to the Proposition 8.4 \cite
{22} there is a correspondence between knot mosaics and grid
diagrams \cite{3,7,24}, that extends to mirror-curves.

The mosaic number $m(L)$ of a link $L$ is the smallest number $n$ for which $L$ is representable as a link $n$-mosaic \cite{22}.

\begin{theorem}
For every link $L$, the mosaic number $m(L)= p+q$, where $p$
and $q$ are dimensions of the minimal $RG[p,q]$ in which $L$ can be
realized. The dimension of the grid (arc) representation
equals $m(L)+1=p+q+1$.
\end{theorem}

\begin{figure}[th]
\centering
{ \includegraphics[scale=0.4]{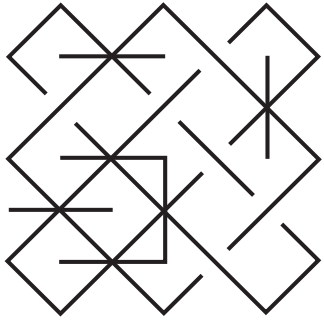}}
\caption{Mirror-curve diagram of the knot $2\,1\,1\,2$ ($6_3$) in
$RG[3,3]$. \label{f1.14}}
\end{figure}

Conjecture 10.4 \cite{22} is an easy corollary of this theorem,
claiming that the mosaic number of the knot $2\,1\,1\,2$ ($6_3$) is
$6$, since its minimal rectangular grid is $RG[3,3]$, and its code is
$\{\{2, -2, 1\}, \{1, 1, -2\}, \{-2, -2, -2\}, \{1, 1, 1\}\}$ (Figure~\ref{f1.14}).

Notice that the knot $2\,1\,1\,2$ ($6_3$) does not cover
$RG[3,3]$ entirely-- if a square in our grid contains just a curl (kink) which can be undone with the Reidemeister I move,  we call it empty square or a hole. Hence, it may be useful to look at the minimal size of every mirror-curve, i.e., the minimal number of non-empty squares necessary to draw it in some (hollow) polyomino \cite{17}.

\begin{conjecture}
Mosaic number of a connected sum $L_1\#L_2$ of two links $L_1$ and $L_2$ satisfies the following equality:
$$m(L_1\#L_2)= m(L_1)+m(L_2)-3.$$
\end{conjecture}

There are two additional numbers that potentially describe the structure of mirror-curves related to the unknotting (unlinking) number:
\begin{itemize}
  \item the minimal number of two-sided mirrors that we
need to add to some mirror-curve in order to obtain unlink,
  \item maximal number of mirrors that can be added to it without obtaining
unlink.
\end{itemize}
 For example, for a $RG[p,2]$ ($p\ge 2$) the first number
equals $p-1$, and the other equals $3p-4$.

\section{Product of mirror-curves}

Algebraic operation called {\it product} can be defined for mirror-curves derived from the same
rectangular grid $RG[p,q]$ by promoting symbols  $2$, $-2$, $1$, and $-1$ in their codes to elements of a semigroup of order $4$ \cite{26}. For example, consider the semigroup $S$ of order $4$, generated by elements $A=\{a,aba\}$,
$B=\{b,bab\}$, $C=\{ab\}$, and $D=\{ba\}$, with the semigroup operation given in the Cayley table:

\begin{center}
\begin{tabular}{|c|c|c|c|c|}
  \hline
  * & $A$ & $B$ & $C$ & $D$ \\ \hline
  $A$ & $A$ & $C$ & $C$ & $A$ \\ \hline
  $B$ & $D$ & $B$ & $B$ & $D$ \\ \hline
  $C$ & $A$ & $C$ & $C$ & $A$ \\ \hline
  $D$ & $D$ & $B$ & $B$ & $D$ \\ \hline

\end{tabular}
\end{center}

\begin{figure}[th]
\centering
{ \includegraphics[scale=0.5]{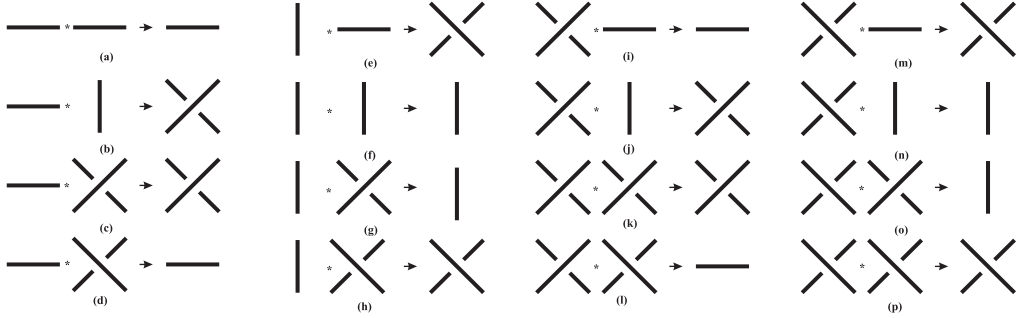}}
\caption{(a) $2*2\rightarrow 2$; (b) $2*-2\rightarrow 1$; (c)
$2*1\rightarrow 1$; (d) $2*-1\rightarrow 2$; (e) $-2*2\rightarrow
-1$; (f) $-2*-2\rightarrow -2$; (g) $-2*1\rightarrow -2$; (h)
$-2*-1\rightarrow -1$; (i) $1*2\rightarrow 2$; (j) $1*-2\rightarrow
1$; (k) $1*1\rightarrow 1$; (l) $1*-1\rightarrow 2$; (m)
$-1*2\rightarrow -1$; (n) $-1*-2\rightarrow -2$; (o)
$-1*1\rightarrow -2$; (p) $-1*-1\rightarrow -1$. \label{f1.15}}
\end{figure}

\begin{figure}[th]
\centering
{ \includegraphics[scale=0.35]{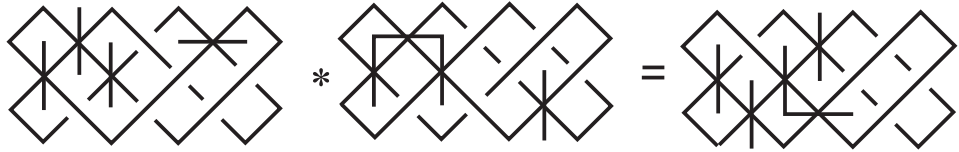}}
\caption{Product
$M_1*M_2=\{\{-2,-2,1,1\},\{2,1\},\{-2,2\},\{-1,-1\}\}$ of
mirror-curves $M_1=\{\{-2,-2,1,1\},$ $\{1,2\},$ $\{-1,1\},$ $\{-1,-2\}\}$ and
$M_2=\{\{-2,-2,1,1\},\{-1,-2\},\{1,-1\},\{2,-1\}\}$. \label{f1.16}}
\end{figure}

First, we substitute $2\rightarrow a$, $-2\rightarrow b$,
$1\rightarrow ab$, $-1\rightarrow ba$, use the semigroup product and then substitute the original
symbols back (Figure~\ref{f1.15}), to obtain the code $M_1*M_2=\{\{-2,-2,1,1\},$ $\{2,1\},\{-2,$ $2\},\{-1,-1\}\}$ as
the product of mirror-curves
$M_1=\{\{-2,-2,1,1\},\{1,2\},\{-1,1\},\{-1,$ $-2\}\}$ and
$M_2=\{\{-2,-2,1,$ $1\},\{-1,-2\},\{1,-1\},\{2,-1\}\}$  (Figure~\ref{f1.16}).

\begin{figure}[th]
\centering
{ \includegraphics[scale=0.35]{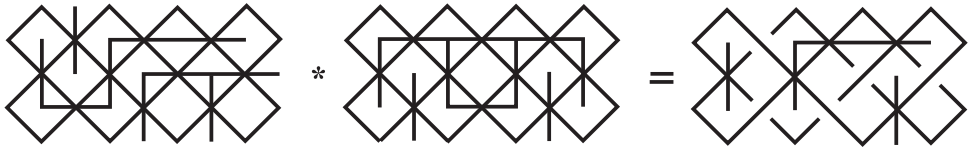}}
\caption{Product $M_1*M_2=\{\{-2,-2,1,1\},\{-1,1\},\{1,-2\},\{2,-2\}\}$ of
mirror-curves $M_1=\{\{-2,-2,2,2\},$ $\{-2,2\},$ $\{2,-2\},$ $\{2,-2\}\}$ and
$M_2=\{\{-2,-2,-2,-2\},$ $\{2,-2\},\{-2,-2\},\{2,-2\}\}$.
\label{f1.17}}
\end{figure}

Since the elements $a$, $b$, $ab$ and $ba$ are idempotents, we have
the equality $M*M=M^2=M$  for every mirror-curve $M$. If
$M_{[p,q]}$ is the set of all mirror-curves derived from $RG[p,q]$,
the basis (minimal set of mirror-curves from which $M_{[p,q]}$ can
be obtained by the operation of product) is the subset of all
mirror-curves of dimensions $p\times q$ with codes consisting only of $2$'s and $-2$'s (Figure~\ref{f1.17}),
i.e. the set of all unlinks belonging to $RG[p,q]$. The basis is not closed under the operation of product: the product of two mirror-curves belonging does not belong to the same basis, since it has at least one crossing.

\begin{figure}[th]
\centering
{ \includegraphics[scale=0.4]{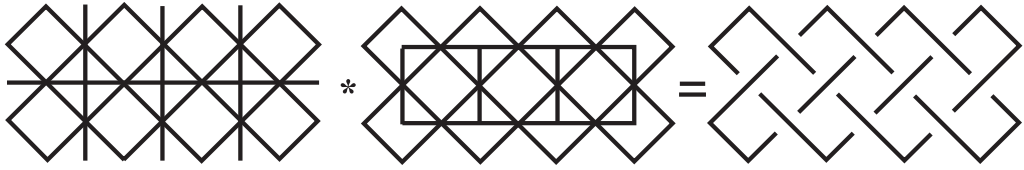}}
\caption{Alternating link $3\,1\,2\,1\,3$ ($L10a_{101}$)
corresponding to $RG[4,2]$ obtained as the product
$M_1*M_2=\{\{1,1,1,1\},\{1,1\},\{1,1\},\{1,1\}\}$ of mirror-curves
$M_1=\{\{2,2,2,2\},\{2,2\},\{2,2\},\{2,2\}\}$ and
$M_2=\{\{-2,-2,-2,-2\},$ $ \{-2,-2\},\{-2,-2\},\{-2,-2\}\}$.
\label{f1.18}}
\end{figure}

In particular, alternating knot or link corresponding to $RG[p,q]$
is obtained as the product of mirror-curves containing
only vertical and horizontal mirrors, see Figure~\ref{f1.18}. Substituting
with elements of different semigroups of order $4$ listed in
\cite{11}, we could obtain different multiplication laws for
mirror-curves.

\section{Kauffman bracket polynomial and mirror-curves}

Let $L$ be any unoriented link diagram. Define the {\it Kauffman state} $S$ of $L$ to be a choice of smoothing for each crossing of $L$ \cite{18,19,20}. There are two choices of smoothing for each crossing, $A$-smoothing and $B$-smoothing, and thus there are $2^c$ states of a diagram with $c$ crossings. In a similar way, we can define the Kauffman state of $RG[p,q]$ as a mirror-curve in $RG[p,q]$ whose code contains only $2$'s and $-2$'s.

Let us consider the set $M^*_{[p,q]}$, called the Kauffman states of $RG[p,q]$, which contains $2^v$ elements corresponding to the choice of mirrors $2$ or $-2$ in the mid-points of $v=2pq-p-q$ internal edges of $RG[p,q]$. Every element of $M^*_{[p,q]}$ can be characterized by the dimensions $p$ and $q$ of
the grid $RG[p,q]$, and another integer $m$ ($0\le m\le 2^v-1$). In order to obtain the matrix code of some mirror-curve  from $(p,q,m)$ code, substitute $0$ by $2$ and $1$ by $-2$ in the binary expansion of $m$ then subdivide the list into $q-1$ lists of length $p$ and $p-1$ lists of length $q$. This code naturally extends to products of mirror-curves. Every mirror-curve $M$ in $RG[p,q]$ can be represented as a product $M=M_1*M_2$ of two mirror-curves $M_1$ and $M_2$ from the set $M^*_{[p,q]}$, hence it can be denoted by a four-number code $(p,q,m,n)$, compounded from codes  $(p,q,m)$ and $(p,q,n)$ of mirror-curves (Kauffman states) $M_1$ and $M_2$, respectively.

For example, the mirror-curve $M$ corresponding to a trefoil knot in $RG[2,2]$ can be represented by the code $(2,2,1,15)$. By expressing numbers $m=1$ and $n=15$ in $4$-digit binary codes, we obtain $\{0,0,0,1\}$ and $\{1,1,1,1\}$, so $M$ is the product of the mirror-curves $\{\{2,2\},\{2,-2\}\}$ and $\{\{-2,-2\},\{-2,-2\}\}$.
Four-number code is not unique. For example, a trefoil in $RG[2,2]$ can be represented by $(2,2,1,15)$, $(2,2,2,15)$, $(2,2,4,15)$, and $(2,2,8,15)$. We choose the minimal code $(2,2,1,15)$ as the code of the trefoil knot.

This approach provides an easy algorithm for computing the Kauffman bracket
polynomial \cite{18,19,20} of an alternating link $L$ directly from its
mirror-curve representation. The Kauffman state sum approach
bypasses the recursive skein relation definition of the Kauffman bracket polynomial, which is given by the formula
$$\sum_S a^{A(S)}a^{-B(S)}(-a^2-a^{-2})^{|S|-1},$$
\noindent as the sum over all Kauffman states $S$ of a link $L$, where $A(S)$ and $B(S)$ is the number of $A$-smoothings and $B$-smoothings, respectively, and $|S|$ is the number of components in the particular state \cite{18,19}.

Analogously, the Kauffman bracket polynomial can be computed as the sum of all possible states of the mirror-curve representing our link $L$.

\begin{figure}[th]
\centering
{ \includegraphics[scale=0.5]{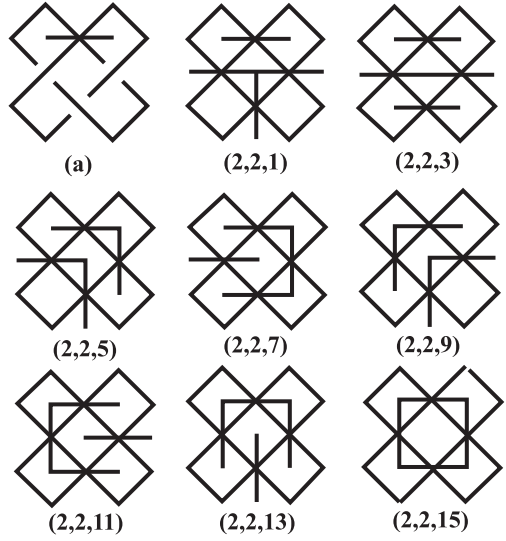}}
\caption{Computation of the Kauffman bracket polynomial for a
trefoil. \label{f1.19}}
\end{figure}

Since all Kauffman states of a link $L$ represented by a mirror-curve $M$ in a grid $R[p,q]$ form a subset of $M^*_{[p,q]}$, the Kauffman bracket polynomial can be computed from the data associated to the mirror-curves in $M^*_{[p,q]}$. Let $M_i$ be a mirror-curve corresponding to some Kauffman state $S_i$ of a link $L$. Denote by $A_i$ be the number of mirrors labeled $1$ in $M$ that changed to $2$ in $M_i$, and $|M_i|$ be the number of components of a Kauffman state $M_i$. Then the bracket polynomial of $L$ can be expressed as
\begin{equation}
  <M>=\sum _{i=0}^{2^n-1} a^{A_i}a^{-n+A_i}(-a^2-a^{-2})^{|M_i|-1}
\end{equation}

For example, a trefoil given by the mirror-curve
$(2,2,1,15)=\{\{1,1\},\{1,-2\}\}$,  shown in Figure~\ref{f1.19}a, has $8$ states:
$\{\{2,2\},\{2,-2\}\}$, $\{\{2,2\},\{-2,-2\}\}$,
$\{\{2,-2\},\{2,-2\}\}$, $\{\{2,-2\},\{-2,$ $-2\}\}$,
$\{\{-2,2\},\{2,-2\}\}$, $\{\{-2,2\},\{-2,-2\}\}$, $\{\{-2,-2\},$ $
\{2,-2\}\}$, $\{\{-2,-2\},$ $\{-2,-2\}\}$ given by the codes
$(2,2,2k+1)$, $0\le k\le 7$, see Figure~\ref{f1.19}.

According to the multiplication table shown on Figure~\ref{f1.15}), a mirror image of a link $L$ given as a product mirror-curve $M=M_1*M_2$,
is $M'=M_2*M_1.$ If $M=M_1*M_2$,  the pair of mirror-curves $(M_1,M_2)$ will be called the {\it
decomposition} of $M$. Minimal decomposition yields the minimal mirror-curve code $(p,q,m,n)$
for every link $L$. For example, the Hopf link is given by the minimal
$(p,q,m,n)$-code $(2,2,1,14)$, trefoil by $(2,2,1,15)$, figure-eight knot by
$(3,2,7,127)$, {\it etc.}

To facilitate computations of the Kauffman bracket polynomial we use two special Kauffman states with all smoothings of one kind: $A$-state ($B$-state) that contain only $A$-smoothings ($B$-smoothings)\footnote{In the language of the Kauffman states of mirror-curves, this means that the first contains only $2$'s, and the other $-2$'s.}.

Let us denote by $M_0=(p,q,0)$ the $A$-state, and by $M_{2^v-1}=(p,q,2^v-1)$ the $B$-state of $RG[p,q].$
\begin{theorem}
Every representation of an alternating link $L$  as a mirror-curve in $RG[p,q]$
can be given as a (left or right) product of some Kauffman state $M$  with $M_0$ or $M_{2^v-1}$, determined by a code $(p,q,m,2^v-1)$ or $(p,q,0,n)$, with $v=2pq-p-q$
and $m,n\in \{0,2^v-1\}$.
\end{theorem}

Such a representation of an alternating link $L$ will be called {\it
canonical representation}. For example, the minimal representation
of the Hopf link is $(2,2,1,14)$, and its canonical representation
is $(2,2,5,15)$. The minimal and canonical
representation of an alternating link cannot always be obtained from the
same rectangular grid. Similarly, the minimal representation of
the knot $4\,2$ ($6_1$) can be obtained from $RG[3,2]$, and its
first canonical representation from $RG[4,2]$.

Every non-alternating mirror-curve $M$ in $RG[p,q]$ can be uniquely represented as the
product of two alternating mirror-curves $M_1=(p,q,m_1,n_1)$ and
$M_2=(p,q,m_2,n_2)$. This means that every non-alternating link $L$ or an alternating link given by its
non-alternating mirror-curve diagram can be denoted by the minimal code of the
form $(p,q,m_1,n_1,m_2,n_2)$.

In order to compute the Kauffman bracket polynomial of  non-alternating
links from mirror-curves we can use the preceding results obtained
for alternating mirror-curves and extend our computation to all
mirror-curves by using skein relation for bracket polynomial, i.e.,
the product of mirror-curves. For example,  consider a
non-alternating link $2,2,-2$ ($6_3^3$) in $RG[3,3]$, given by the
code $M=\{\{1, 1, -1\}, \{1, 1, -1\}, \{-2, 2, -2\}, \{-2, 2,
-2\}\}$.
Let  $<M>$ denote the bracket polynomial of the mirror-curve $M$. Then:
\begin{eqnarray*}
  <M>&=&a(a<M_0>+a^{-1}<M_1>)+a^{-1}(a<M_2>+a^{-1}<M_3>)\\
     &=& a^2<M_0>+<M_1>+<M_2>+a^{-2}<M_3>,
\end{eqnarray*}
\noindent where
\begin{eqnarray*}
  M_0&=&\{\{1, 1, -2\}, \{1, 1, -2\}, \{-2, 2, -2\}, \{-2, 2, -2\}\},\\
  M_1&=&\{\{1, 1, -2\}, \{1, 1, 2\}, \{-2, 2, -2\}, \{-2, 2, -2\}\},\\
  M_2&=&\{\{1, 1, 2\}, \{1, 1, -2\}, \{-2, 2, -2\}, \{-2, 2, -2\}\},\\
  M_3&=&\{\{1, 1, 2\}, \{1, 1, 2\}, \{-2, 2, -2\},\{-2, 2, -2\}\}
\end{eqnarray*}
\noindent are mirror-curves with all crossings positive. Hence,
\begin{eqnarray*}
  <M>&=&a^2(2+a^{-8}+a^8)+(-a^{-6}-a^2+a^6-a^{10})+(-a^{-6}-a^2+a^6-a^{10})+\\
  && a^{-2}(1+a^{-8}+a^{-4}+a^{12}) = a^{-10}+a^{-2}+2a^6.
\end{eqnarray*}

Notice that we have used all Kauffman states, this time expanded over all negative crossings. In the
case of a non-alternating mirror-curve $M$ with $n$ crossings, and $n_-$ negative
crossings the Kauffman bracket polynomial is given by the following state sum formula:
\begin{equation}
  <M>=\sum _{i=0}^{2^n-1}a^{A_i}a^{-n_-+A_i}<M_i>,
\end{equation}
\noindent where $A_i$ is the number of mirrors changed from $1$ in $M$ to
$-2$ in a Kauffman state $M_i$, and $M_i$ ($0\le i\le 2^n-1$) are
alternating mirror-curves obtained as the Kauffman states taken over
negative crossings by changing $-1$ into $-2$ and $2$. Since every mirror-curve $M_i$  corresponding to some Kauffman state $S_i$ is just a collection of $|M_i|=|S_i|$ circles, its Kauffman bracket is $<M_i>=(-a^2-a^{-2})^{|S_i|-1}.$  Moreover, the power of $a^{A_i}a^{-n_-+A_i}$ is the vertex weight $w_i$: the number of $A$-smoothings minus the number of $B$ smoothings in a state $S_i$ times $\pm 1$, depending on the sign of each crossing.
 The state sum formula for the Kauffman bracket polynomial \cite{19} now has the following form:
\begin{equation}
  \sum _{i=0}^{2^n-1}a^{w_i}(-a^2-a^{-2})^{|S_i|-1}
\end{equation}

\begin{example}
 To illustrate the formula above, we give an explicit computation of the Kauffman bracket using the formula above, for the mirror-curve $M=\{\{1,1\},\{-1,-2\}\}$ shown on Figure~\ref{f1.20}a, which is just an unknot represented as a trefoil with one crossing change.

\begin{figure}[th]
\centering
{ \includegraphics[scale=0.4]{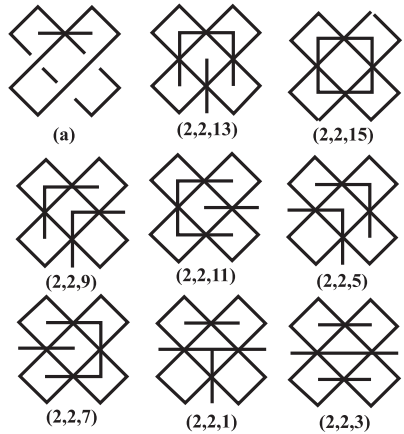}}
\caption{Computation of the Kauffman bracket polynomial for the
mirror-curve $M=\{\{1,1\},\{-1,-2\}\}$ (a) and its eight states.
\label{f1.20}}
\end{figure}

Eight mirror-curves $M_i$ corresponding to the Kauffman states $S_i$, $i=0,\ldots ,7$ are shown
on Figure~\ref{f1.20} and their codes, as well as the number of components, are contained in Table~\ref{tri}.
Next we compute the vertex weights $(w_0,\ldots, w_7)=(3,1,1,-1,1,-1,-1,-3),$  to obtain $<M>=-a^3$.

\begin{table}
\begin{center}
  \begin{tabular}{|l |c|}
  \hline
  Kauffman state $M_i=S_i$ & $|M_i|=|S_i|$\\  \hline
  $M_0=\{\{-2,-2\},\{2,-2\}\}$ & 1 \\   \hline
  $M_1=\{\{-2,-2\},\{-2,-2\}\}$ & 2 \\  \hline
   $M_2=\{\{-2,2\},\{2,-2\}\}$&  2\\  \hline
 $M_3=\{\{-2,2\},\{-2,-2\}\}$  &  1\\  \hline
   $M_4=\{\{2,-2\},\{2,-2\}\}$ &  2\\  \hline
   $M_5=\{\{2,-2\},\{-2,-2\}\}$&  1\\  \hline
 $M_6=\{\{2,2\},\{2,-2\}\}$  &  3\\  \hline
  $M_7=\{\{2,2\},\{-2,-2\}\}$ & 2 \\ \hline
\end{tabular}
\end{center}
\caption{ $2^n$ mirror-curves $M_i$ ($i=0,\ldots ,2^n-1$) shown on Figure~\ref{f1.20} and the number of their link components.}
\label{tri}
\end{table}
\end{example}

\section{L-polynomials and mirror-curves}

Mirror-curves can also be used for computing the Kauffman L-polynomial
\cite{19,20} defined by the following axioms:

\begin{enumerate}
\item $L(+1)+ L(-1)= z(L(0)+L(\infty));$
\item $L \looparrowright =aL;$
\item $L \looparrowleft =a^{-1}L;$
\item $L(\bigcirc )=1;$
\end{enumerate}

\noindent where  $\looparrowright $ and $\looparrowleft $
denote positive and negative curls.

\begin{figure}[th]
\centering
{ \includegraphics[scale=0.8]{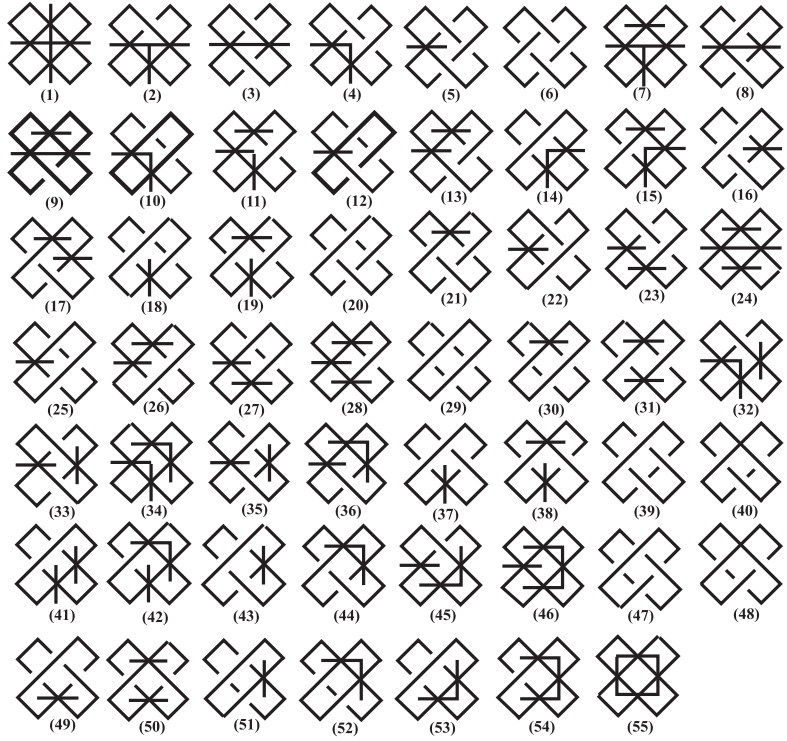}}
\caption{Mirror-curves obtained from $RG[2,2]$ up to isometry.
\label{f1.21}}
\end{figure}

Grid $RG[2,2]$ contains $55$ mirror-curves\footnote{Up to isometry.}
shown on Figure~\ref{f1.21}), where the mirror-curves (20) and (47) reduce to (44), (23) and
(52) reduce to (55), (50) reduces to (24), and (30) reduces to (54).
Knowing that $L(\bigcirc ^n)=\delta ^{n-1}$, where $\delta
=({{{a+a^{-1}} \over z}}-1)$, we can compute the L-polynomial
for all of them except for the mirror-curves (6), (21), (31) and (44)
 by simply counting circles and curls.

\begin{figure}[th]
\centering
{ \includegraphics[scale=0.55]{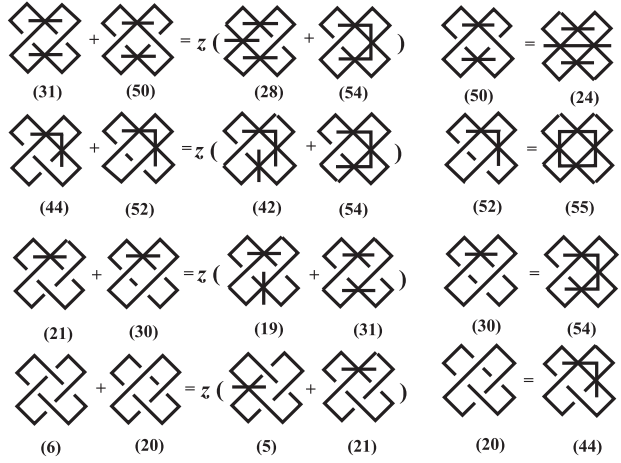}}
\caption{Computation of L-polynomial from mirror-curves in
$RG[2,2]$. \label{f1.22}}
\end{figure}

We have the following relations which are also illustrated on Figure \ref{f1.22}:
\begin{eqnarray*}
&&L(\{\{1,1\},\{-2,-2\}\})+L(\{\{1,-1\},\{-2,-2\}\})=z(L(\{\{1,2\},\{-2,-2\}\})+L(\{\{1,-2\},\{-2,-2\}\})),
\end{eqnarray*}
\noindent with $L(\{\{1,-1\},\{-2,-2\}\})=L(\{\{2,2\},\{-2,-2\}\}).$

In other words, we have $L(31)+L(50)=z(L(28)+L(54))$, with $L(50)=L(24)$.
\begin{eqnarray*}
&&L(\{\{1,-2\},\{1,-2\}\})+L(\{\{1,-2\},\{-1,-2\}\})=z(L(\{\{1,-2\},\{2,-2\}\})+L(\{\{1,-2\},\{-2,-2\}\})),
\end{eqnarray*}
\noindent where $L(\{\{1,-2\},\{-1,-2\}\})=L(\{\{-2,-2\},\{-2,-2\}\}),$ i.e., $L(44)+L(52)=z(L(42)+L(54))$, with $L(52)=L(55)$:
\begin{eqnarray*}
&&L(\{\{1,1\},\{1,-2\}\})+L(\{\{1,1\},\{-1,-2\}\})=z(L(\{\{1,1\},\{2,-2\}\})+L(\{\{1,1\},\{-2,-2\}\})).
\end{eqnarray*}
\noindent with $L(\{\{1,1\},\{-1,-2\}\})=L(\{\{1,-2\},\{-2,-2\}\}),$ i.e., $L(21)+L(30)=z(L(19)+L(31))$. Since $L(30)=L(54)$
 we have
\begin{eqnarray*}
  &&L(\{\{1,1\},\{1,1\}\})+L(\{\{1,1\},\{1,-1\}\})=z(L(\{\{1,1\},\{1,2\}\})+L(\{\{1,1\},\{1,-2\}\})),
\end{eqnarray*}
\noindent with $L(\{\{1,1\},\{1,-1\}\})=L(\{\{1,-2\},\{1,-2\}\})$ i.e., $L(6)+L(20)=z(L(5)+L(21))$, with $L(20)=L(44)$, see Figure~\ref{f1.22}.

\noindent Hence, we conclude that
\begin{eqnarray*}
  L(Hopf\,Link)&=&L(2_1^2)= L(31)=z(L(28)+L(54))-L(24)=z(a^{-1}+a)-\delta ^2\\
               &=&-(a^{-1}+a)z^{-1}+1+(a^{-1}+a)z,\\
  L(Right\,Trefoil )&=&L(3_1)= L(21)=z(L(19)+L(31))-L(54)=z(a^{-2}+L(31))-a=  \\
            &=& -(a^{-1}+2a)+(a^{-2}+1)z+(a^{-1}+a)z^2,\\
   L(4_1^2) &=&L(6)=z(L(5)+L(21))-L(44)=z(a^{-3}+L(21))-L(44) \\ &=&-(a^{-1}+a)z^{-1}-1+(a^{-3}-2a^{-1}-3a)z+(a^{-2}+1)z^2+(a^{-1}+a)z^3.
\end{eqnarray*}

\begin{figure}[th]
\centering
{ \includegraphics[scale=0.5]{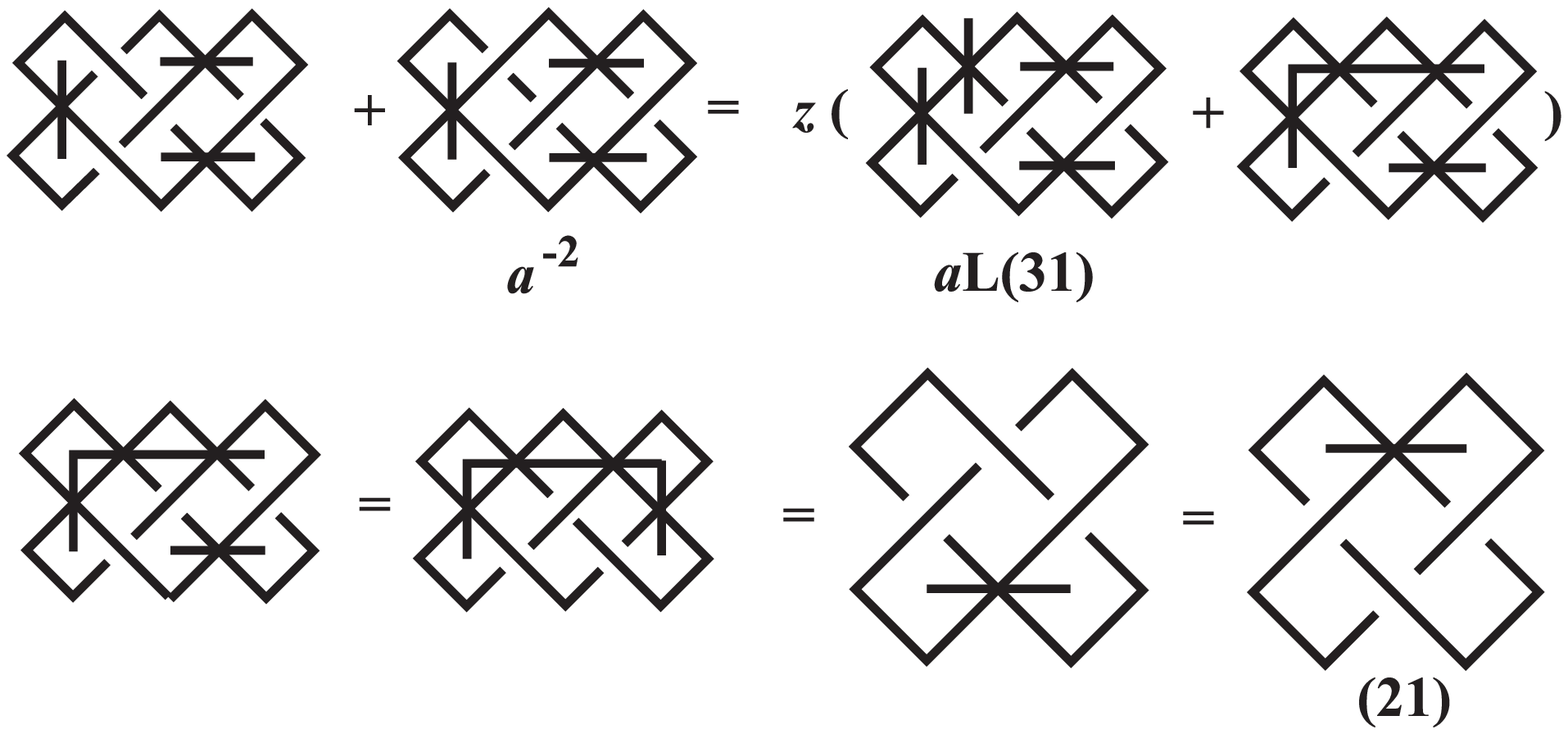}}
\caption{Computation of L-polynomial for figure-eight knot.
\label{f1.23}}
\end{figure}

In general, L-polynomials for mirror-curves can be computed in the same way, or by simplifying
computations using previously obtained results and relations.
 For example, the L-polynomial of the mirror-curve, see Figure~\ref{f1.23},
$\{\{-2,1,1\},\{1,1\},\{-2,-2\}\}$ which represents the figure-eight
knot $4_1$ in $RG[3,2]$  satisfies the relation:
\begin{eqnarray*}
&&L(\{\{-2,1,1\},\{1,1\},\{-2,-2\}\})+L(\{\{-2,1,1\},\{1,-1\},\{-2,-2\}\}))=\\
&&z(L(\{\{-2,1,1\},\{1,2\},\{-2,-2\}\})+L(\{\{-2,1,1\},\{1,-2\},\{-2,-2\}\}))
\end{eqnarray*}
Since $L(\{\{-2,1,1\},\{1,-1\},\{-2,-2\}\})=a^{-2}$,
$L(\{\{-2,1,1\},\{1,2\},\{-2,$ $-2\}\})=aL(31)$, and the
mirror-curve $\{\{-2,1,1\},\{1,-2\},\{-2,-2\}\}$ reduces to
$\{\{1,1\},\{-2,1\}\}=\{\{1,1\},\{-2,1\}\}$, i.e., to the
mirror-curve (21) in $RG[2,2]$ corresponding to the trefoil knot,
\begin{eqnarray*}
   L(4_1)&=&z(aL(31)+L(21))-a^{-2}=(-a^{-2}-1-a^2)-(a^{-1}+a)z+\\
         & & (a^{-2}+2+a^2)z^2+(a^{-1}+a)z^3.
\end{eqnarray*}

This approach can also be used for deriving recursive formulas relating
the L-polynomials of knot and link families given in Conway
notation. Members of the knot family $p$ ($p\ge 1$), denoted by
Conway symbols as $1$, $2$, $3$, $4$, $5$, $\ldots $, namely the unknot, Hopf link $2_1^2$, trefoil $3_1$, link $4_1^2$, knot $5_1$,
$\ldots $  satisfy the following recursion:
\begin{eqnarray*}
  L(1)&=&a\\
  L(2)&=&-(a^{-1}+2a)+(a^{-2}+1)z+(a^{-1}+a)z^2 \\
  L(p)&=&z(a^{-p+1}+L(p-1))-L(p-2), \,\,\, \text { for }  p\ge 3.\\
\end{eqnarray*}

For the knot family $p\,2$ ($p\ge 2$), which consists from knots
$4_1$, $5_2$, $6_1$, $7_2$, $\ldots $ we have the recursion
\begin{eqnarray*}
 L(1\,2)&=& L(3)\\
  L(2\,2)&=&(-a^{-2}-1-a^2)-(a^{-1}+a)z+(a^{-2}+2+a^2)z^2+(a^{-1}+a)z^3\\
  L(p\,2)    &=&z(L((p-1)\,2)+a^{p-1}L(2))-L((p-2)\,2), \,\,\, \text{ for } p\ge 3.
\end{eqnarray*}
Members of the link family $3\,p$ ($p\ge 3$) satisfy the recursion

$$L(3\,p)=z(L(2\,p)+a^2L(p))-L(p+1), \,\,\, \text{ for } p\ge 3,$$

\noindent where $2\,p$ is the mirror image of the link $p\,2$, and
$3\,p$ is the mirror image of the link $p\,3$.

In general, the link family $p\,q$ ($p\ge q\ge 2$) satisfies the following
recursion

$$L(p\,q)=z(L((p-1)\,q)+a^{p-1}L(q))-L((p-2)\,q).$$

\noindent {\bf Acknowledgement:} The authors express their gratitude to
the Ministry of Science and Technological Development for providing
partial support for this project (Grant No. 174012).

\section{Appendix}
\begin{center}

\begin{table}\tiny
\begin{center}
   \begin{tabular}{|c|c|c|c|}   \hline
 1  & $3\,1\,2\,1\,3$ & $L10a_{101}$ & $\{\{-1,-1,-1,-1\},\{-1,-1\},\{-1,-1\},\{-1,-1\}\}$ \\  \hline
 2  & $5\,1\,3$ & $9_5$ &  $\{\{-1,2,-1,-1\},\{-1,-1\},\{-1,-1\},\{-1,-1\}\}$ \\ \hline
 3  & $3\,1\,2\,1\,2$ & $9_{20}$ &  $\{\{-2,-1,-1,-1\},\{-1,-1\},\{-1,-1\},\{-1,-1\}\}$ \\ \hline
 4  & $4\,1\,1\,3$ & $9_5^2$ & $\{\{-1,-1,1,1\},\{-1,-1\},\{-1,-1\},\{1,1\}\}$ \\ \hline
 5  & $3\,1\,3\,2$ & $9_8^2$ & $\{\{-1,-1,-1,1\},\{-1,-1\},\{-1,-1\},\{1,1\}\}$ \\  \hline
 6  & $3\,1\,1\,1\,3$ & $9_9^2$ & $\{\{-1,-1,-1,-1\},\{-1,-1\},\{-2,-1\},\{-1,-1\}\}$ \\  \hline
 7  & $5\,1\,2$ & $8_2$ & $\{\{-2,-1,2,-1\},\{-1,-1\},\{-1,-1\},\{-1,-1\}\}$ \\  \hline
 8  & $4\,1\,3$ & $8_4$ & $\{\{-2,2,-1,-1\},\{-1,-1\},\{-1,-1\},\{-1,-1\}\}$ \\  \hline
 9  & $3\,1\,1\,1\,2$ & $8_{13}$ & $\{\{-2,-1,-1,-1\},\{-1,-1\},\{-2,-1\},\{-1,-1\}\}$ \\  \hline
 10  & $8$ & $8_1^2$ & $\{\{-1,2,2,-1\},\{-1,-1\},\{-1,-1\},\{-1,-1\}\}$ \\  \hline
 11  & $4\,2\,2$ & $8_3^2$ &  $\{\{-1,1,-1,-1\},\{-1,-1\},\{1,1\},\{-1,-1\}\}$ \\  \hline
 12  & $3\,2\,3$ & $8_4^2$ & $\{\{-1,-1,-1,-1\},\{-1,-1\},\{-1,1\},\{-1,-1\}\}$ \\  \hline
 13  & $3\,1\,2\,2$ & $8_5^2$ & $\{\{-2,-1,-1,-1\},\{-2,-1\},\{-1,-1\},\{-1,-1\}\}$ \\  \hline
 14  & $2\,4\,2$ & $8_6^2$ & $\{\{-1,1,1,-1\},\{-1,-1\},\{1,1\},\{-1,-1\}\}$ \\  \hline
 15  & $2\,1\,2\,1\,2$ & $8_7^2$ &  $\{\{-2,-1,-1,-2\},\{-1,-1\},\{-1,-1\},\{-1,-1\}\}$ \\  \hline
 16  & $7$ & $7_1$ & $\{\{-2,2,2,-1\},\{-1,-1\},\{-1,-1\},\{-1,-1\}\}$  \\  \hline
 17  & $5\,2$ & $7_2$ & $\{\{-2,-1,2,-1\},\{-2,-1\},\{-1,-1\},\{-1,-1\}\}$ \\  \hline
 18  & $2\,2\,1\,2$ & $7_6$ & $\{\{-2,-1,-1,-2\},\{-2,-1\},\{-1,-1\},\{-1,-1\}\}$ \\  \hline
 19  & $2\,1\,1\,1\,2$ & $7_7$ & $\{\{-2,-1,-1,-2\},\{-1,-1\},\{-2,-1\},\{-1,-1\}\}$ \\  \hline
 20  & $4\,1\,2$ & $7_1^2$ & $\{\{-2,2,-1,-2\},\{-1,-1\},\{-1,-1\},\{-1,-1\}\}$ \\  \hline
 21  & $3\,1\,1\,2$ & $7_2^2$ & $\{\{-1,-1,-1,-2\},\{-1,-1\},\{-2,-1\},\{-2,-1\}\}$ \\  \hline
 22  & $2\,3\,2$ & $7_3^2$ &  $\{\{1,-1,-1,1\},\{-1,-1\},\{-1,-1\},\{1,1\}\}$ \\  \hline
 23  & $2\,1\,1\,2$ & $6_3$ & $\{\{-2,-1,-1,-2\},\{-1,-1\},\{-2,-1\},\{-2,-1\}\}$ \\  \hline
 24  & $3\,3$ & $6_2^2$ & $\{\{-1,1,-1,-1\},\{-1,-1\},\{-1,-1\},\{-1,-1\}\}$ \\  \hline
 25  & $2\,2\,2$ & $6_3^2$ & $\{\{-2,-1,-1,-2\},\{-2,-1\},\{-1,-1\},\{-2,-1\}\}$ \\  \hline
\end{tabular}
\end{center}
\caption{ $KL$s derived from $RG[4,2]$}
\end{table}
\begin{figure}[th]
\centering
{ \includegraphics[scale=0.7]{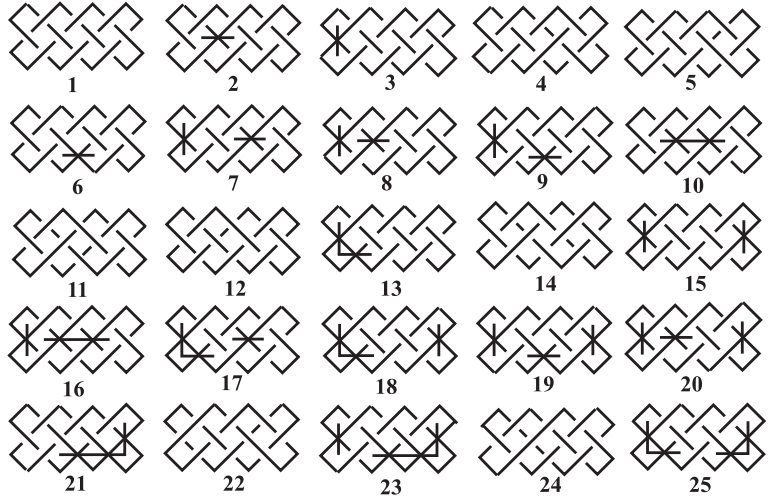}}
\caption{Mirror-curves $1$-$25$ derived from $RG[4,2]$. \label{f1.24}}
\end{figure}

\begin{table}\begin{center}

\tiny
  \begin{tabular}{|c|c|c|c|}   \hline

1  & $8^*2:2:2:2$ &  &
$\{\{-1,-1,-1\},\{-1,-1,-1\},\{-1,-1,-1\},\{-1,-1,-1\}\}$ \\  \hline

2 & $(2,2)\,(3\,1,-3\,1)$ &  &
$\{\{-1,1,1\},\{-1,-1,-1\},\{-1,-1,-1\},\{1,1,-1\}\}$ \\  \hline

3  & $(-5\,1,2)\,(2,2)$ &  &
$\{\{-1,-1,1\},\{-1,-1,-1\},\{-1,-1,-1\},\{1,-1,-1\}\}$ \\  \hline

4  & $6^*-2.2.-2:4$  &  &
$\{\{-1,-1,1\},\{1,-1,-1\},\{-1,-1,1\},\{1,-1,-1\}\}$ \\  \hline

5  & $6^*3.2.-3:2$ &  &
$\{\{1,1,1\},\{-1,-1,-1\},\{1,-1,-1\},\{1,-1,-1\}\}$ \\  \hline

6  & $6^*-3.-3\, 0::-3\, 0$ &  &
$\{\{-1,1,-1\},\{-1,1,-1\},\{-1,1,-1\},\{-1,1,-1\}\}$ \\  \hline

7  & $2\,1\,2\,1\,1\,1\,2$ & $10_{44}$ &
$\{\{-1,-1,-1\},\{-2,-1,-1\},\{-1,-1,-1\},\{-1,-2,-1\}\}$ \\  \hline

8  & $.4.2\,0$ & $10_{85}$ &
$\{\{-1,-1,-1\},\{-2,-1,-1\},\{-1,-1,-1\},\{-1,2,-1\}\}$ \\  \hline

9  & $4\,1\,2\,1\,2$ & $L10a_{99}$ &
$\{\{-1,-1,-1\},\{-1,-2,-1\},\{-1,-1,-1\},\{-1,2,-1\}\}$ \\  \hline

10  & $.3:3\,0$ & $L10a_{140}$ &
$\{\{1,1,1\},\{-1,-1,-1\},\{1,1,-1\},\{1,-1,-1\}\}$ \\  \hline

11  & $6,2,2$ & $L10a_{145}$ &
$\{\{-1,-1,-1\},\{-1,2,-1\},\{-1,-1,-1\},\{-1,2,-1\}\}$ \\  \hline

12  & $.2.3.2\,0$ & $L10a_{162}$ &
$\{\{-1,-1,1\},\{-1,-1,-1\},\{-1,-1,-1\},\{1,1,-1\}\}$ \\  \hline

13  & $8^*2::2$ & $L10a_{163}$ &
$\{\{-1,-1,-1\},\{-2,-1,-1\},\{-1,-1,-1\},\{-2,-1,-1\}\}$ \\  \hline

14  & $2\,0.2.2\,0.2\,0$ & $L10a_{164}$ &
$\{\{-1,-1,1\},\{-1,-1,-1\},\{-1,-1,-1\},\{-1,-1,-1\}\}$ \\  \hline

15  & $(2\,1,-2\,1)\,(2,2)$ & $L10n_{73}$ &
$\{\{1,1,1\},\{1,-1,-1\},\{1,-1,-1\},\{1,-1,-1\}\}$ \\  \hline

16  & $(3\,1,-2)\, (2,2)$ & $L10n_{85}$ &
$\{\{-1,1,1\},\{-1,-1,-1\},\{-1,-1,-1\},\{-1,1,-1\}\}$ \\  \hline

17  & $(2,2)\,(4,-2)$ & $L10n_{86}$ &
$\{\{-1,-1,1\},\{-1,-1,-1\},\{1,-1,-1\},\{1,-1,-1\}\}$ \\  \hline

18  & $4,3\,1,-2$ & $L10n_{92}$ &
$\{\{-1,-1,1\},\{1,-1,-1\},\{-1,-1,-1\},\{1,-1,-1\}\}$ \\  \hline

19  & $4,4,-2$ & $L10n_{93}$ &
$\{\{-1,1,-1\},\{-1,-1,-1\},\{-1,1,-1\},\{-1,1,-1\}\}$ \\  \hline

20  & $2\,0.-2.-2\,0.2\,0$ & $L10n_{94}$ &
$\{\{-1,-1,1\},\{1,-1,-1\},\{1,-1,1\},\{1,-1,-1\}\}$ \\  \hline

21  & $3\,1,3\,1,-2$ & $L10n_{95}$ &
$\{\{-1,-1,-1\},\{-1,-1,-1\},\{-1,1,-1\},\{-1,-1,-1\}\}$ \\  \hline

22  & $4\,1\,2\,2$ & $9_{11}$ &
$\{\{-1,-1,-1\},\{-2,-2,-1\},\{-1,-1,-1\},\{-1,2,-1\}\}$ \\  \hline

23  & $4\,1\,1\,1\,2$ & $9_{14}$ &
$\{\{-2,-1,-1\},\{-1,-2,-1\},\{-1,-1,-1\},\{-1,2,-1\}\}$ \\  \hline

24  & $2\,1\,3\,1\,2$ & $9_{17}$ &
$\{\{-1,-1,-1\},\{-2,-1,-1\},\{-1,2,-1\},\{-1,-2,-1\}\}$ \\  \hline

25  & $2\,2\,1\,2\,2$ & $9_{23}$ &
$\{\{-1,-1,-1\},\{-2,-2,-1\},\{-1,-1,-1\},\{-1,-1,-2\}\}$ \\  \hline

26  & $2\,1\,2\,1\,1\,2$ & $9_{27}$ &
$\{\{-1,-1,-1\},\{-2,-1,-2\},\{-1,-1,-1\},\{-1,-2,-1\}\}$ \\  \hline

27  & $2\,1\,1\,1\,1\,1\,2$ & $9_{31}$ &
$\{\{-2,-1,-1\},\{-1,-2,-1\},\{-1,-1,-1\},\{-2,-1,-1\}\}$ \\  \hline

28  & $6\,1\, 2$ & $9_1^2$ &
$\{\{-2,-1,-1\},\{-1,2,-1\},\{-1,-1,-1\},\{-1,2,-1\}\}$ \\  \hline

29  & $2\,2\,1\,1\,1\,2$ & $9_{12}^2$ &
$\{\{-2,-1,-1\},\{-1,-1,-2\},\{-1,-1,-1\},\{-1,-2,-1\}\}$ \\  \hline

30  & $5,2,2$ & $9_{13}^2$ &
$\{\{-1,-1,-1\},\{-2,2,-1\},\{-1,-1,-1\},\{-1,2,-1\}\}$ \\  \hline

31  & $.4$ & $9_{31}^2$ &
$\{\{-2,-1,-1\},\{-1,2,-1\},\{-1,-1,-1\},\{-2,-1,-1\}\}$ \\  \hline

32  & $.3.2\,0$ & $9_{35}^2$ &
$\{\{-2,-1,-1\},\{-1,-1,-2\},\{-1,-1,-1\},\{-1,2,-1\}\}$ \\  \hline

33  & $8^*2$ & $9_{42}^2$ &
$\{\{-2,-1,-1\},\{-1,-1,-2\},\{-1,-1,-1\},\{-2,-1,-1\}\}$ \\  \hline

34  & $6\,2$ & $8_1$ &
$\{\{-1,-1,-1\},\{-1,2,-1\},\{-1,2,-1\},\{-2,-2,-1\}\}$ \\  \hline

35  & $3,3,2$ & $8_5$ &
$\{\{-1,-1,-1\},\{-2,-1,-1\},\{-1,2,-1\},\{-2,2,-1\}\}$ \\  \hline

36  & $4\,1\,1\,2$ & $8_7$ &
$\{\{-2,-1,-1\},\{-2,-2,-1\},\{-1,-1,-1\},\{-1,2,-1\}\}$ \\  \hline

37  & $2\,3\,1\,2$ & $8_8$ &
$\{\{-2,-1,-1\},\{-1,-1,-2\},\{-1,2,-1\},\{-1,-2,-1\}\}$ \\  \hline

38  & $2\,1,3,2$ & $8_{10}$ &
$\{\{-2,-1,-1\},\{-1,-1,-2\},\{-2,-1,-1\},\{-1,2,-1\}\}$ \\  \hline

39  & $2\,2\,2\,2$ & $8_{12}$ &
$\{\{-2,-1,-1\},\{-1,-2,-1\},\{-2,-1,-1\},\{-1,-1,-2\}\}$ \\  \hline

40  & $2\,2\,1\,1\,2$ & $8_{14}$ &
$\{\{-2,-2,-1\},\{-1,-1,-2\},\{-1,-1,-1\},\{-2,-1,-1\}\}$ \\  \hline

41  & $.2.2\,0$ & $8_{16}$ &
$\{\{-2,-1,-1\},\{-2,2,-1\},\{-1,-1,-1\},\{-1,-1,-2\}\}$ \\  \hline

42  & $.2.2$ & $8_{17}$ &
$\{\{-2,-1,-1\},\{-2,-1,-2\},\{-1,-1,-1\},\{-1,-1,-2\}\}$ \\  \hline

43  & $8^*$ & $8_{18}$ &
$\{\{-2,-1,-1\},\{-2,-1,-2\},\{-1,-1,-1\},\{-2,-1,-1\}\}$ \\  \hline

44  & $2\,1\,2\,1\,2$ & $8_7^2$ &
$\{\{-2,-1,2\},\{-1,-1,-2\},\{-1,-1,-1\},\{-2,-1,-1\}\}$ \\  \hline

45  & $2\,1\,1\,1\,1\,2$ & $8_8^2$ &
$\{\{-2,-1,-1\},\{-2,-2,-1\},\{-1,-1,-1\},\{-2,-1,-1\}\}$ \\  \hline

46  & $4,2,2$ & $8_1^3$ &
$\{\{-2,2,-1\},\{-1,-1,-2\},\{-1,-1,-1\},\{-1,2,-1\}\}$ \\  \hline

47  & $3\,1,2,2$ & $8_2^3$ &
$\{\{-1,-1,1\},\{-1,-1,-1\},\{1,-1,-1\},\{-1,-1,-1\}\}$ \\  \hline

48  & $(2,2)\,(2,2)$ & $8_4^3$ &
$\{\{-2,-1,-1\},\{-1,-1,-2\},\{-2,-1,-1\},\{-1,-1,-2\}\}$ \\  \hline

49  & $.3$ & $8_5^3$ &
$\{\{-2,-1,-1\},\{-2,-1,-2\},\{-1,-1,-1\},\{-1,2,-1\}\}$ \\  \hline

50  & $.2:2\,0$ & $8_6^3$ &
$\{\{-1,-1,-1\},\{-2,-1,-2\},\{-1,-1,-1\},\{-2,-1,-2\}\}$ \\  \hline

51  & $4,2,-2$ & $8_7^3$ &
$\{\{-1,-1,1\},\{1,-1,-1\},\{-1,1,-1\},\{1,-1,-1\}\}$ \\  \hline

52  & $3\,1,2,-2$ & $8_8^3$ &
$\{\{-1,-1,1\},\{-1,-1,-1\},\{-1,1,-1\},\{-1,-1,-1\}\}$ \\  \hline

53  & $(2,2)\,(2,-2)$ & $8_9^3$ &
$\{\{1,-1,1\},\{1,-1,-1\},\{1,-1,-1\},\{-1,-1,-1\}\}$ \\  \hline

54  & $(2,2)\,-(2,2)$ & $8_{10}^3$ &
$\{\{-1,-1,-1\},\{-1,-1,-1\},\{-1,1,-1\},\{-1,1,-1\}\}$ \\  \hline

55  & $4\,3$ & $7_3$ &
$\{\{-2,-1,-1\},\{-2,-2,-1\},\{-2,-1,-1\},\{-1,2,-1\}\}$ \\  \hline

56  & $3\,2\,2$ & $7_5$ &
$\{\{-2,-2,-1\},\{-1,-1,-2\},\{-1,-1,-1\},\{-2,-1,-2\}\}$ \\  \hline

57  & $2,2,2+$ & $7_1^3$ &
$\{\{-1,1,1\},\{-1,-1,-1\},\{1,1,-1\},\{-1,1,-1\}\}$ \\  \hline

58  & $2\,3\,2$ & $7_3^2$ &
$\{\{-2,-2,-1\},\{-1,2,-2\},\{-1,-1,-1\},\{-2,-1,-1\}\}$ \\  \hline

59  & $3,2,2$ & $7_4^2$ &
$\{\{-2,-1,-1\},\{-1,-1,-2\},\{-1,2,-1\},\{-2,2,-1\}\}$ \\  \hline

60  & $2\,1,2,2$ & $7_5^2$ &
$\{\{-2,-1,-1\},\{-2,-1,-2\},\{-2,-1,-1\},\{-1,-1,-2\}\}$ \\  \hline

61  & $.2$ & $7_6^2$ &
$\{\{-2,-1,-2\},\{-2,2,-1\},\{-1,-1,-1\},\{-1,-1,-2\}\}$ \\  \hline

62  & $2,2,2$ & $6_1^3$ &
$\{\{-2,-1,-1\},\{-2,-1,-2\},\{-2,-1,-1\},\{-2,-1,-2\}\}$ \\  \hline

63  & $6^*$ & $6_2^3$ &
$\{\{-2,-1,-1\},\{-2,-1,-2\},\{-1,-1,-1\},\{-2,2,-2\}\}$ \\  \hline

64  & $2,2,-2$ & $6_3^3$ &
$\{\{-1,-1,1\},\{-1,1,-1\},\{1,-1,-1\},\{-1,-1,-1\}\}$ \\  \hline

\end{tabular}\end{center}
\normalsize
\caption{$KL$s derived from $RG[3,3]$}
\end{table}

\begin{figure}[th]
\centering
\subfigure
{ \includegraphics[scale=0.54]{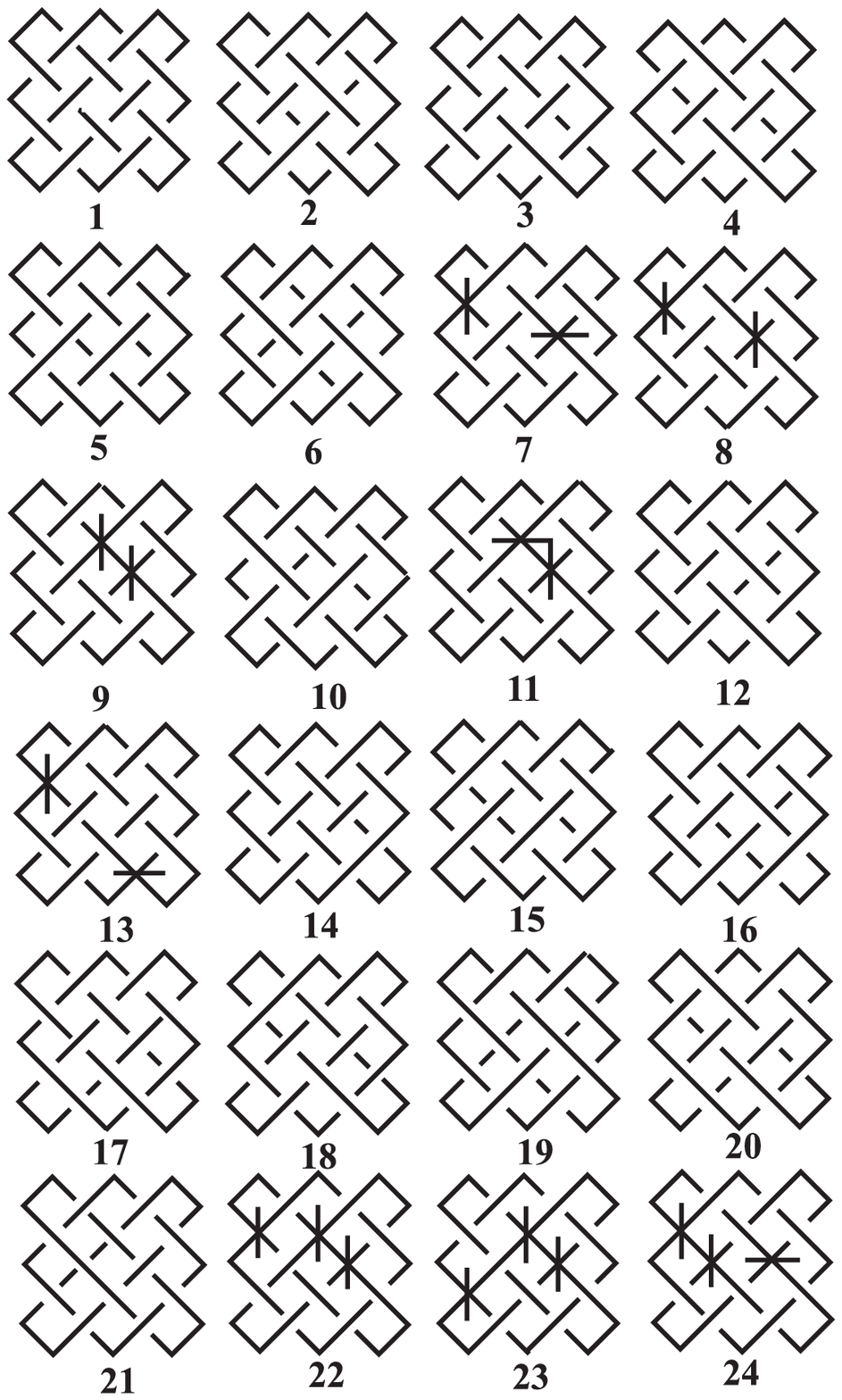}}
\subfigure
{\includegraphics[scale=0.6]{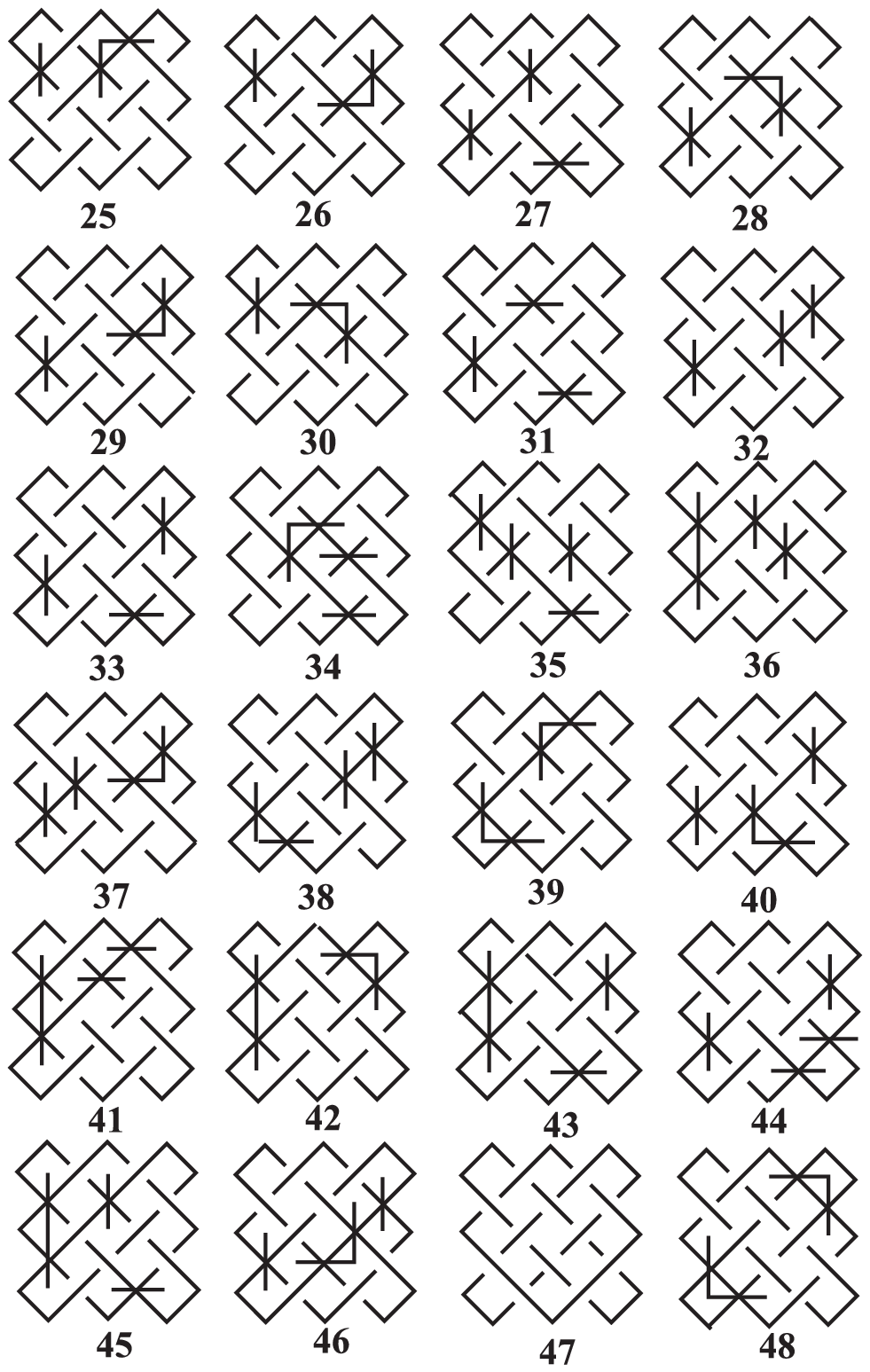} }
\subfigure{
\includegraphics[scale=0.5]{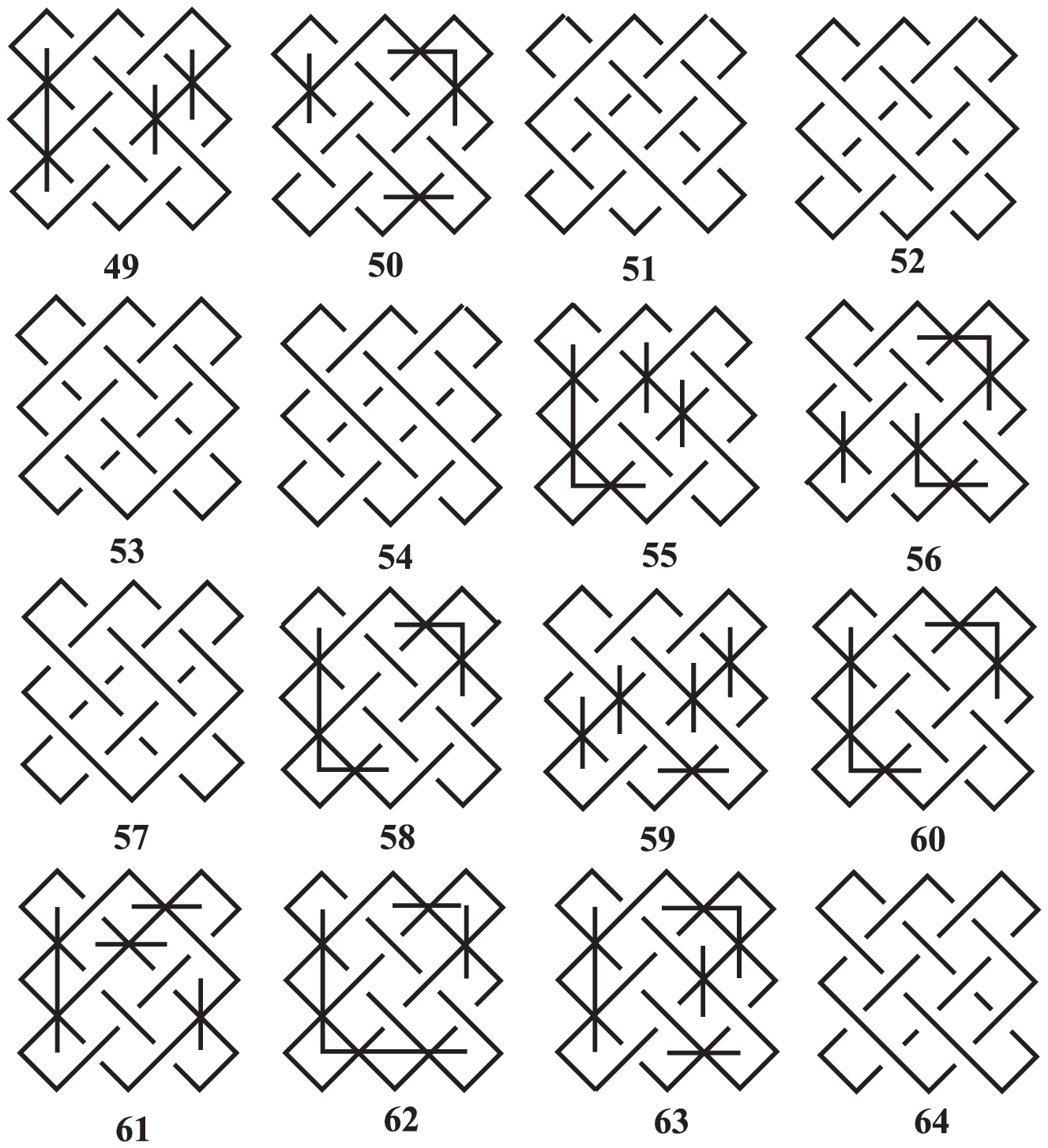}}
\caption{Mirror-curves $1$-$64$ derived from $RG[3,3]$. \label{f1.27}}
\end{figure}
\end{center}

\end{document}